\documentstyle[12pt,psfig]{article}

\def\interior{{\rm Int}}
\def\vbar{{\overline v}}

\def\Kbar{{\overline K}}
\def\GL{{\rm GL}}
\def\defin#1#2{\smallskip\noindent{\bf Definition (#1).} #2}
\def\Ccell{{\rm C}^{\rm cell}}
\def\Cphi{{\rm C}^\varphi}
\def\Xw{X_{\rm w}}

\def\tauw{\tau_{\rm w}}

\def\basis{{\hbox{\Got b}}}
\def\hbasis{{\hbox{\Got h}}}
\def\orient{{\hbox{\Got o}}}

\def\move#1{\mathop{
\begin{picture}(10,5)
\put(0,-3.8){$\tilde{ }$}
\put(1,1){$\to$}
\end{picture}}
\limits^{#1}}

\def\stwotriv{S^2_{\rm triv}}
\def\bthtr{B^3_{\rm triv}}

\reversemarginpar

\makeatletter
\def\@begintheorem#1#2{\it \trivlist \item[\hskip \labelsep{\bf #1\ #2.}]}
\newtheorem{teo}{Theorem}[section]
\newtheorem{rem}[teo]{Remark}
\newtheorem{lem}[teo]{Lemma}

\newtheorem{prop}[teo]{Proposition}

\def\finedim#1{{\hfill\hbox{\enspace\fbox{\ref{#1}}}}\vspace{5pt}}
\def\dim#1{\vspace{1pt}\noindent{\it Proof of} {\hspace{2pt}}\ref{#1}.}

\def\cont{{\rm C}}
\def\combd{{\rm Comb}^\partial}
\def\objd{{\rm Obj}^\partial}
\def\recd{r^\partial}
\def\leg{{\rm Leg}}

\font\scpicc=cmcsc10
\font\sc=cmcsc10 scaled 1200

\newfont{\Bbb}{msbm10 scaled 1200}
\def\mr{{\hbox{\Bbb R}}}

\def\mz{{\hbox{\Bbb Z}}}

\newfont{\mycal}{eusm10 scaled 1200}
\def\ttt{{\hbox{\mycal T}}}

\newfont{\Got}{eufm10 scaled 1200}

\def\vv{{\cal V}}

\def\tt{{\cal T}}

\def\ee{{\cal E}}
\def\rr{{\cal R}}

\def\pp{{\cal P}}

\def\dd{{\cal D}}

\def\trasvint{\cap\hspace{-8.5pt}|\hspace{5pt}}

\title{Torsion Invariants of Combed 3-Manifolds with Boundary Pattern
and Legendrian Links}

\author{Riccardo Benedetti\qquad Carlo Petronio\thanks{The second named
author gratefully acknowledges financial support by GNSAGA-CNR}}

\begin{document}

\maketitle

\noindent{\small{\scpicc Abstract}. We extend Turaev's definition of torsion
invariants of 3-dimensional manifolds equipped with 
non-singular vector fields, by allowing (suitable)
tangency circles to the boundary, and manifolds with non-zero Euler
characteristic. We show that these invariants apply in particular to (the
exterior of) Legendrian links in contact 3-manifolds. Our approach uses a
combinatorial encoding of vector fields, based on standard spines. In this
paper we extend this encoding from closed manifolds to manifolds with
boundary.}

\vspace{.5cm}


\vspace{.5cm}

\noindent Reidemeister torsion is a classical yet very vital topic in
3-dimensional topology, and it was recently used in a variety of important
developments. To mention a few, torsion is a fundamental ingredient of the
Casson-Walker-Lescop invariants (see {\em e.g.}~\cite{lescop}), and more
generally of the perturbative approach to quantum invariants (see {\em
e.g.}~\cite{lemuoh}). Relations have been pointed out between torsion and
hyperbolic geometry~\cite{porti}. Turaev's torsion of non-singular vector
fields on 3-manifolds~\cite{turaev:Euler} has been recognized to have deep
connections with some 3-dimensional versions of the  Seiberg-Witten
invariants~\cite{meng:taub},~\cite{turaev:spinc}. It is also worth recalling
that vector fields (and framings), have also been used by Kuperberg~\cite{kuperberg} to
construct new invariants of different nature, based on Hopf algebras more
general than those employed for quantum invariants. (There are reasons to
speculate that also Kuperberg's  invariants should have a torsion content
and relations with Turaev's work, but we do not insist on this.)

In this paper, using (actually, improving) our theory of branched standard
spines~\cite{lnm} and building on~\cite{turaev:Euler}, we extend Turaev's
definition of torsion by allowing vector fields to have (appropriate) tangency
circles to the boundary. Moreover, we do not require the manifolds to have zero
Euler characteristic (an assumption which is at the base of Turaev's theory). 
To be precise, we accept any compact oriented manifold with (possibly empty)
boundary, and non-singular vector fields with the only restriction that their
orbits should be tangent to the boundary ``from inside'' ({\it i.e.}~along
circles of {\em concave} apparent contour). Equivalence is given by homotopy
through fields of the same type. Recall that Turaev never accepts tangency to
the boundary, and the most important applications  in~\cite{turaev:spinc} are
given when the boundary is empty or a union of tori with the field pointing
inwards on them.

One of the topological situations in which we are able to define torsion
invariants arises quite naturally when one considers Legendrian links in
contact 3-manifolds. Remark in particular that our torsions are defined also
for homologically non-trivial Legendrian links, when the usual invariants, such
as the rotation number (Maslov index), are not defined.

Our definition of torsion is based on a combinatorial encoding of non-singular
concave vector fields on manifolds with boundary. The topological-geometric
counterpart of this encoding is the theory of standard spines with the further
structure of branched surface. The foundations of this theory, together with
the combinatorial encoding of vector fields in the special case of closed
manifolds, were provided in~\cite{lnm}. The main technical result of the
present paper is the extension of this encoding to the case with boundary. We
also show that one of the combinatorial moves taken into account in~\cite{lnm}
can be discarded. This implies the rather significant fact that the branched
versions of the basic 2-to-3 Matveev-Piergallini move are actually
sufficient.  

If a concave vector field $v$ on a manifold $M$ is encoded by a branched spine
$P$, in order to compute torsions we consider a cell complex $X(P)$, whose support
is  obtained by collapsing the whole $\partial M$ to an internal point. Moreover
we use $v$ to construct a ``canonical spider'' in $X(P)$, and then we use this
spider to determine a fundamental family of cells in the universal cover of
$X(P)$. The fact that we are able to accept manifolds with non-zero Euler
characteristic essentially depends on the fact that our complex naturally arises
as a {\em pointed} space, and the spider is connected, with head at the basepoint.
We recover the situation considered by Turaev as a variation on our basic
definition, by leaving uncollapsed the boundary components on which the field
points inwards. The use of branched standard spines, in connection with torsion,
allows a considerable simplification of the proof of invariance. All our proofs
will be direct combinatorial arguments. Even if obviously inspired to
Turaev's~\cite{turaev:Euler}, our work is essentially self-contained.

We conclude the introduction by pointing out an interesting subtlety which arises
when dealing with torsions (see Subsection~\ref{comparison:warning:subsection} for
exact statements). We start by recalling that, in general, the torsion of a
complex $X$, which depends on a ring homomorphism
$\varphi:\mz[\pi_1(X)]\to\Lambda$, is an element of $\Kbar_1(\Lambda)$ which is
only well-defined up the action of  $\varphi(\pi_1(X))$. Turaev's main achievement
in~\cite{turaev:Euler} is the proof that this action can be disposed of when $X$
is a manifold with a  vector field defined on it. However, it may still be the
case that the automorphism group of $X$ acts non-trivially on torsions. We show in
the present paper that when dealing with Legendrian links this action can
sometimes be neglected, which leads to a sharper version of the invariant.

\section{Torsion(s) of a branched spine}

In this section we will briefly recall the notion of branched spine, we will
describe how to associate a certain cell complex $X(P)$ to each branched spine
$P$, and we will show that the branching of $P$ allows to define a canonical
``spider'' in $X(P)$. We will then define a torsion 
$\tau^{\varphi}(P,\hbasis)\in\Kbar_1(\Lambda)$ for a (suitable) representation
$\varphi$ of $\pi_1(X)$ into the multiplicative group of a (suitable) ring
$\Lambda$, where $\Kbar_1(\Lambda)$ denotes the Whitehead group of $\Lambda$
and $\hbasis$ is a $\Lambda$-basis of the $\varphi$-twisted homology module of
$X(P)$, which is assumed to be free. Later we will describe some variations on
the definitions of $X(P)$ and $\tau^{\varphi}(P,\hbasis)$. 

The first subsection establishes various notations used extensively in this
paper.

\subsection{Reminder on branched spines}\label{reminder}

A quasi-standard polyhedron $P$ is a  finite connected 2-dimensional polyhedron
with singularity of stable nature (triple lines and points where 6 non-singular
components meet). Such a $P$ is called {\it standard} if all the components of
the natural stratification given by singularity are open cells. Depending on
dimension, we will call the components {\it vertices, edges} and {\it regions}.
A {\it screw-orientation} for $P$ is a screw-orientation on all its edges, with
the obvious compatibility at vertices.  A {\it branching} on $P$ is an
orientation for each region in $P$, such that no edge is induced the same
orientation three times. See~\cite{lnm} for careful definitions of all these notions.

A standard polyhedron with a fixed oriented branching and screw-orientation
will be called {\em branched polyhedron}, and denoted typically by $P$. We will not
use specific notations for the extra structures: they will be considered to be
part of $P$. Unless the contrary is explicitly stated, by ``manifold'' we will
mean a connected, oriented, compact 3-dimensional manifold, with or without
boundary. Using the {\it Hauptvermutung}, we will freely intermingle the
differentiable, piecewise linear and topological viewpoints. Homeomorphisms
will always respect orientations.

All vector fields mentioned in this paper will be non-singular, and they will be
termed just {\em fields} for the sake of brevity. On a manifold $M$ we will
consider {\em concave} fields, namely non-singular fields which are simply tangent
to $\partial M$ only along a finite union $\Gamma$ of circles, with the property
that near tangency points orbits are contained in $\interior(M)$. Note that
$\Gamma$ splits $\partial M$ into a {\em white} portion $W$, on which the field
points inwards, and a {\em black} portion $B$, on which the field points outwards.

It turns out~\cite{manuscripta} that  a branched polyhedron with a
screw-orientation is automatically the spine of a manifold, which is unique
by~\cite{casler}, so we will often replace the term `polyhedron' by `spine'.
The following result, proved in~\cite{lnm}, is the starting
point of our constructions.

\begin{prop}\label{from:spine:to:field}
To every branched spine $P$ there corresponds a manifold $M(P)$ 
with non-empty boundary and a concave field $v(P)$ on $M(P)$.
The pair $(M(P),v(P))$ is well-defined up to homeomorphism. Moreover
an embedding $i:P\to\interior(M(P))$ is defined, 
and has the property that $v(P)$ is positively transversal to $i(P)$.
\end{prop}

\begin{figure}
\centerline{\psfig{file=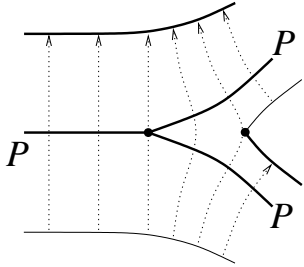,width=4cm}}
\caption{\label{constr:M} Manifold and field associated to a branched spine.}
\end{figure}

The topological construction which underlies this proposition is actually quite
simple, and it is illustrated in Fig.~\ref{constr:M}. Concerning the last
assertion of the proposition, note that the branching allows to define an
oriented tangent plane at each point of $P$.

Even if it is not necessary now, we inform the reader that suitable
restrictions of the map $P\mapsto(M(P),v(P))$  defined above are surjective.
One of the main achievements of~\cite{lnm} was the introduction of an
equivalence relation on branched spines  which makes this map injective, when
restricted to manifolds bounded by $S^2$ (essentially, closed manifolds). In
Section~\ref{improve:lnm:section} we will provide a substantial improvement of this
result.

\paragraph{Spines and ideal triangulations.} We remind the reader that an {\it
ideal triangulation} of a manifold $M$ with non-empty boundary is  a partition
$\tt$ of $\interior(M)$ into open cells of dimensions 1, 2 and 3, induced by a
triangulation $\tt'$ of the space $Q(M)$, where:
\begin{enumerate}
\item $Q(M)$ is obtained from $M$ by collapsing each component of
$\partial M$ to a point;
\item $\tt'$ is a triangulation only in a loose sense, namely
self-adjacencies and multiple adjacencies of tetrahedra are allowed;
\item The vertices of $\tt'$ are precisely the points of $Q(M)$
which correspond to components of $\partial M$.
\end{enumerate}

\begin{figure}
\centerline{\psfig{file=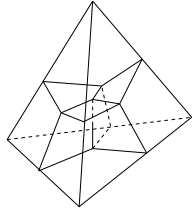,width=2.5cm}}
\caption{\label{duality} Duality between standard spines and ideal
triangulations.}
\end{figure}

It turns out (\cite{mafo},~\cite{tesi},~\cite{matv:new}) that there exists a natural
bijection between standard spines and ideal triangulations
of a 3-manifold. Given an ideal triangulation, the
corresponding standard spine is just the 2-skeleton of the dual
cellularization, as illustrated in Figure~\ref{duality}.
The inverse of this correspondence will be denoted by $P\mapsto\tt(P)$. 
It will be convenient in the sequel to call {\it centre} of a cell of $P$
the only point in which the cell meets the simplex of $\tt(P)$ dual to it.

Now let $P$ be a branched spine. First of all, we can realize $\tt(P)$ in such
a way that its edges are orbits of the restriction of $v(P)$ to
$\interior(M(P))$, and the 2-faces are unions of such orbits.  Being orbits,
the edges of $\tt(P)$ have a natural orientation, and the branching condition
implies (as remarked in~\cite{gr:1}) that on each tetrahedron of $\tt(P)$ 
exactly one of the vertices is a sink and one is a source. 

It is quite interesting to remark that not only the edges, but also the faces
and the tetrahedra of $\tt(P)$ have natural orientations. For tetrahedra, we
just restrict the orientation of $M(P)$. For faces, we first note that the edges
of $P$ have a natural orientation (the prevailing orientation induced by the
incident regions). Now, we orient a face of $\tt(P)$ so that the algebraic
intersection in $M(P)$ with the dual edge is positive. 

\subsection{Basic definition of torsion}\label{basic:subsection}

Let $P$ be a branched spine.

\defin{triangulated complex associated to $P$}{Let $X(P)=M(P)/\partial(M(P))$
and let $x_0(P)$ be the point which corresponds to $\partial(M(P))$. Since
$X(P)\setminus\{x_0(P)\}\cong \interior(M(P))$, $P$ naturally embeds in $X(P)$.
Moreover, a field $\vbar(P)$ is defined on $X(P)\setminus\{x_0(P)\}$, and the
closure in $X(P)$ of each infinite half-orbit of $\vbar(P)$ is obtained by
adding $x_0(P)$. The ideal triangulation $\tt(P)$ of $M(P)$ induces a
triangulation $\ttt(P)$ of $X(P)$, with only vertex $x_0(P)$ and open
positive-dimensional simplices which correspond to those of $\tt(P)$ and are
unions of orbits of $\vbar(P)$.}

\begin{rem}{\em
\begin{enumerate}
\item[(i)] The Euler characteristic of $X(P)$ is given by 
$$\chi(X(P))=1-(1/2)\cdot\chi(\partial(M(P)))=1-\chi(M(P))=1-\chi(P).$$
\item[(ii)] Stipulating $x_0(P)$ to be positive, and using the remarks made above,
we see that also in $\ttt(P)$ all the  
simplices have a natural orientation.
\end{enumerate}
}\end{rem}

\defin{spider associated to $P$}{We define $s(P)$ as the singular 1-chain in
$X(P)$ obtained as $\sum_c\alpha_c$, where $c$ runs over the centres of cells
of $P$, and $\alpha_c$ is the closure of the positive orbit of $\vbar(P)$ which
starts at $c$. Note that the final endpoint of each $\alpha_c$ is $x_0(P)$.}

\begin{rem}\label{non:zero:chi:good:1}{\em Let $\varepsilon(c)=(-1)^d$ if $c$ is
the centre of a $d$-cell. Then
$$\partial\left(\sum\nolimits_c\varepsilon(c)\cdot\alpha_c\right)=
(1-\chi(X(P)))\cdot x_0(P)+\sum\nolimits_c\varepsilon(c)\cdot c.$$ 
So, if $\chi(X(P))=0$,
the chain $\sum_c\varepsilon(c)\alpha_c$ is precisely a spider (or {\em Euler
chain})  in Turaev's terminology~\cite{turaev:Euler}. But our definition makes
sense also for $\chi(X(P))\neq0$, so our situation is indeed more general. 
We have a preferred basepoint which  is part of the
structure from the beginning,  and the ``error'' in the boundary of our spider is
automatically located at the basepoint. This fact will be crucial also
below.}\end{rem}

For the sake of brevity, in the sequel we will denote $\pi_1(X(P),x_0(P))$ just by
$\pi$. We denote now by $(\tilde X(P),\tilde x_0(P))$ a universal cover of
$(X(P),x_0(P))$. The reason for considering pointed spaces is that any two such
covers  are {\em canonically} isomorphic, and all our constructions will obviously
be equivariant under this isomorphism. On $\tilde X(P)$ we consider the action of
$\pi$ defined using the basepoint $\tilde x_0(P)$. We denote by $\tilde\ttt(P)$
the $\pi$-invariant lifting of $\ttt(P)$ to $\tilde X(P)$. We will consider in the
sequel the complex  $\Ccell_*(\tilde X(P);\mz)$ of integer chains in $\tilde
X(P)$  which are cellular with respect to $\tilde\ttt(P)$. In a natural way,
$\Ccell_*(\tilde X(P);\mz)$ is a complex  of $\mz[\pi]$-modules. Moreover each
$\Ccell_i(\tilde X(P);\mz)$ is free, and a free basis is determined by the choice
of an ordering for the $i$-simplices in $\ttt(P)$ and one lifting for each of them
(as remarked, orientations are canonical).

\defin{lifted spider and free basis}{We define $\tilde s(P)$ as the singular
1-chain $\sum_c\tilde\alpha_c$ in $\tilde X(P)$, where  $\tilde\alpha_c$ is the
only lifting of $\alpha_c$ with second endpoint $\tilde x_0(P)$. We choose
$\tilde x_0(P)$ as preferred lifting of $x_0(P)$. For a positive-dimensional
simplex of $\ttt(P)$  dual to a cell with centre $c$, we choose as preferred
lifting the one which contains the first endpoint of $\tilde\alpha_c$. If
$\sigma$ is an ordering of the simplices in $\ttt(P)$, we denote by
$g_i(P,\sigma)$ the free $\mz[\pi]$-basis of $\Ccell_i(\tilde X(P);\mz)$
obtained from $\sigma$ and these preferred liftings.}

We briefly review now the general algebraic machinery used to define
torsions~\cite{milnor}. We consider a ring $\Lambda$ with unit, with the property
that  if $n$ and $m$ are distinct positive integers then $\Lambda^n$ and 
$\Lambda^m$ are not isomorphic as $\Lambda$-modules. We recall that the
Whitehead group $K_1(\Lambda)$ is defined as the
Abelianization of $\GL_\infty(\Lambda)$. Moreover, $\Kbar_1(\Lambda)$ is the
quotient of $K_1(\Lambda)$ under the action of $-1\in\GL_1(\Lambda)=
\Lambda_*\subset\Lambda$. 

Given a free $\Lambda$-module $M$ and two finite bases $\basis=(b_k) $,
$\basis'=(b'_k)$ of $M$, the assumption on $\Lambda$ guarantees that
$\basis $ and $\basis '$ have the same number of elements, so there
exists an invertible square matrix $(\lambda^h_k)$ such that
$b'_k=\sum_\beta\lambda^h_k b_h$. We will denote by
$[\basis '/\basis ]$ the image of $(\lambda^h_k)$ in
$K_1(\Lambda)$.

\defin{twisted homology and chain basis}{Going back to the topological
situation, let us consider now a  group homomorphism $\varphi:\pi\to\Lambda_*$,
and its natural extension $\tilde\varphi:\mz[\pi]\to\Lambda$ (a ring
homomorphism). We can define now the twisted chain complex $\Cphi_*(P)$, where
$\Cphi_i(P)$ is defined as  $\Lambda\otimes_{\tilde\varphi}\Ccell_i(\tilde
X(P);\mz)$, and the boundary operator is induced from the ordinary boundary.
Note that $\Cphi_i(P)$ is a free $\Lambda$-module, and
each $\mz[\pi]$-basis of $\Ccell_i(\tilde X(P);\mz)$ determines a
$\Lambda$-basis of $\Cphi_i(P)$. We will denote by $g^\varphi_i(P,\sigma)$ the
$\Lambda$-basis of $\Cphi_i(P)$ corresponding to $g_i(P,\sigma)$,
and by $H^\varphi_i(P)$ the $i$-th homology group of the complex $\Cphi_*(P)$.
The canonical isomorphism which exists between
two pointed universal covers of $(X(P),x_0(P))$ induces an isomorphism of
the corresponding homology groups, so $H_*^\varphi(P)$ is intrinsically
defined.}

Our assumptions on $\Lambda$ easily imply the following:

\begin{lem}\label{if:acyclic:then:chi}
If $H^\varphi_i(P)=0$ for all $i$ then $\chi(P)=0$.
\end{lem}

\defin{torsion --- acyclic case}{Assume that $H_*^\varphi(P)=0$. Then we can
apply the general definition of torsion of an acyclic chain complex of free
$\Lambda$-modules with assigned bases. We briefly review this definition,
confining ourselves to the case where the boundary modules are free (in
general, a stable basis should be used). So, let $\basis _i$ be a finite subset
of $\Cphi_i(P)$ such that $\partial\basis _i$ is a basis of
$\partial\Cphi_i(P)$. The complex being acyclic, $(\partial\basis
_{i+1})\cdot\basis _i$ is now a basis of $\Cphi_i(P)$, so we can compare it
with $g_i^\varphi(P,\sigma)$. We define
$$\tau_0^\varphi(P,\sigma)=\prod_{i=0}^3\Big[ (\partial\basis
_{i+1})\cdot\basis _i \;\Big/\;g_i^\varphi(P,\sigma)\Big]^{(-1)^i}\in
K_1(\Lambda).$$ (Independence of the $\basis_i$'s and invariance under
isomorphism of pointed universal covers is readily checked.) Of course $\sigma$
is responsible of at most a sign change, so
$\tau^\varphi(P)=\pm\tau_0^\varphi(P,\sigma)\in\Kbar_1(\Lambda)$ is
well-defined.}

\defin{torsion --- general case}  It follows from Lemma~\ref{if:acyclic:then:chi}
that $\Cphi_*(P)$ is often not acyclic. It is a general fact that torsion can be
defined also in this case, provided the homology $\Lambda$-modules are free and
have assigned bases. Namely, if $\hbasis_i$ is a $\Lambda$-basis of
$H^\varphi_i(P)$, we replace $(\partial\basis_{i+1})\cdot\basis_i$ in the above
formula by $(\partial\basis_{i+1})\cdot\tilde\hbasis_i\cdot\basis_i$, where
$\tilde\hbasis_i$ is a lifting of $\hbasis_i$ to $\Cphi_i(P)$.
So, we have a  torsion
$\tau^\varphi(P,\hbasis)\in\Kbar_1(\Lambda)$ when the  $\Lambda$-modules
$H^\varphi_*(P)$ are free with basis $\hbasis$.

It is maybe appropriate here to remark that the choice of a basis $\hbasis$  of
$H^\varphi_*(P)$ and the definition of $\tau^\varphi(P,\hbasis)$ implicitly assume
a description of the universal cover of $X(P)$, which is typically undoable in
practical cases. However,
if one starts from a representation of $\pi$ into the units of a {\em commutative}
ring  $\Lambda$, one can use from the very beginning the maximal Abelian rather
than the universal cover, which makes computations more feasible.

\defin{sign-refined torsion}
An enhancement of the definition of torsion, due to Turaev, applies in our
situation. As remarked above, the sign
ambiguity in the definition of torsion is only due to the ordering $\sigma$ of
the cells of $X(P)$. This ambiguity can be removed by considering a homological
orientation $\orient$ of $X(P)$, namely an orientation of all the spaces 
$H_i(X(P);\mr)$. We briefly recall how this construction goes. Using $\sigma$,
we get bases of the $\Ccell_i(X(P);\mr)$'s. Now we choose bases of the
$H_i(X(P);\mr)$'s compatible with $\orient$, and we compute the torsion of
$\Ccell_*(X(P);\mr)$ using these bases, thus getting $a\in K_1(\mr)=\mr_*$. It
is easily checked that ${\rm sgn}(a)\cdot \tau_0^\varphi(P,\sigma) \in
K_1(\Lambda)$ is now independent of $\sigma$. Hence we get a torsion
$\tau_0^\varphi(P,\hbasis,\orient)\in K_1(\Lambda)$. Of course, in the acyclic
case, $\hbasis$ is omitted.

\paragraph{Manifolds with white boundary components (Turaev's torsion).} It
could happen that in $\partial(M(P))$ there are some white boundary components,
{\it i.e.}~components on which the field $v(P)$ points inwards. In this
case we can modify the definition of the triangulated complex associated to $P$
by identifying together only the boundary components which are not white. We
denote by $\Xw(P)$ the space thus obtained, by $x_{{\rm w},0}(P)$ the basepoint
obtained by collapsing the non-white boundary components, and by $W(P)$ the 
(homeomorphic) image in $\Xw(P)$ of the union of white boundary components. A
field $\vbar_{\rm w}(P)$ is naturally induced by $v(P)$ on 
$\Xw(P)\setminus(W(P)\cup\{x_{{\rm w},0}(P)\})$, and a spider $s_{\rm w}(P)$ with
head in $x_{{\rm w},0}(P)$ can be defined exactly as above.

Now, if we consider the universal covering $\tilde\Xw(P)$ with a basepoint
$\tilde x_{{\rm w},0}$, we can still lift the spider by locating its head in 
$\tilde x_{{\rm w},0}$, but now the legs of the spider only determine liftings
of the cells of $\Xw(P)\setminus W(P)$. So the spider defines a basis of the
relative chain complex $\Ccell_*(\tilde\Xw(P),\tilde W(P);\mz)$ as a
$\mz[\pi_1(\Xw(P))]$-module. Now the construction proceeds exactly as above,
and allows to define a torsion $\tauw^\varphi(P)\in\Kbar_1(\Lambda)$ for
every  homomorphism $\varphi:\pi_1(\Xw(P))\to\Lambda_*$ such that the twisted
chain complex $\Cphi_*(\Xw(P),W(P))$ is acyclic.
Following the scheme mentioned above, one can also define refined torsions
$\tau_{{\rm w},0}^\varphi(P,\hbasis,\orient)\in K_1(\Lambda)$.

An interesting special case of the definition of $\tau_{{\rm w},0}$ is when
$\partial(M(P))$ consists of a (possibly empty) union of white tori together
with one sphere with  black-white splitting consisting of two discs. In this
case $\Xw(P)$ is a combed manifold $N$, closed or bounded by white tori.
Moreover  our spider defines a combinatorial Euler structure in Turaev's
sense~\cite{turaev:Euler} (see Remark~\ref{non:zero:chi:good:1}), and one sees
that torsions coincide. We will see below in Section~\ref{improve:lnm:section}
that every combed manifold bounded by white tori arises as $\Xw(P)$ for some
$P$, so we actually do cover all situations considered by Turaev. See
Subsection~\ref{discussion} for further comments.

\subsection{How to use torsions to compare
objects}\label{comparison:warning:subsection} A theoretical problem  arises when
one wants to use torsions to distinguish objects. We carefully describe this
problem in our setting, but the phenomenon is general.

We start with the following remark. Let $P_0$ and $P_1$ be branched standard
spines, and let $f:P_0\to P_1$ be a homeomorphism which preserves all the
structures of branched spine. For $i=0,1$, consider the combed manifold 
$(M(P_i),v(P_i))$ determined by $P_i$. As already mentioned, this
manifold is only determined up to homeomorphism, but we fix a definite
representative, together with an embedding of $P_i$ into $M(P_i)$
satisfying the usual conditions. Now $f$ extends to a homeomorphism 
$M(P_0)\to M(P_1)$ well-defined up isotopy, which induces
a homeomorphism $F:X(P_0)\to X(P_1)$ again determined up to isotopy.
Of course we have
$$\tau_0^\varphi(P_0,\hbasis,\orient)=
\tau_0^{F_*(\varphi)}(P_1,F_*(\hbasis),F_*(\orient))$$
for all torsions $\tau_0^\varphi(P_0,\hbasis,\orient)$. However, it may in
principle happen for some other homeomorphism $g:M(P_0)\to M(P_1)$ inducing
$G:X(P_0)\to X(P_1)$ that 
$$\tau_0^\varphi(P_0,\hbasis,\orient)\neq
\tau_0^{G_*(\varphi)}(P_1,G_*(\hbasis),G_*(\orient)).$$

To describe the situation in a different way, suppose that we have branched
spines $P_0$ and $P_1$, we already know that $M(P_0)$ and $M(P_1)$ are
homeomorphic, and we want to use torsions to decide whether $P_0$ and $P_1$ are
isomorphic or not. (Of course, this is not a good example, since the
isomorphism problem for branched spines is easy, while the homeomorphism
problem for manifolds is hard, but more appropriate examples will arise later.)
The above remark implies that torsions cannot be used directly, because an
action of the automorphism (actually, mapping class) group of $M(P_0)=M(P_1)$
has to be taken into account. Therefore, we have an analogue of the
Teichm\"uller {\it vs.} moduli space situation: the basic definition of torsion
involves a ``marking'' of the manifold, so, to get a marking-independent
torsion, one has to quotient out under an action of the mapping class group. We
will privilege in this paper this moduli-type approach, but we will mention
in Remark~\ref{one:manifold} and Subsection~\ref{Teich:version:subsection}
situations where a marking is natural, so a Teichm\"uller-type approach is
feasible.

\section{Invariance under sliding}\label{sliding:invariance:section}

The fundamental move for standard spines, which allows to obtain from each
other any two spines of the same manifold, is the Matveev-Piergallini
move~\cite{matv:mossa},~\cite{piergallini}, which corresponds, in terms of the
\begin{figure}
\centerline{\psfig{file=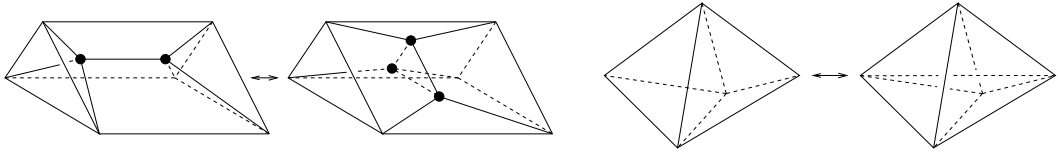,width=14cm}}
\caption{\label{mossa}The fundamental move for spines.}
\end{figure}
dual ideal triangulation, to the so-called 2-to-3 move. Both versions of the
move are illustrated in Fig.~\ref{mossa}. We will consider to be {\em positive}
the move in the direction which increases the number of vertices (or
tetrahedra).

\begin{figure}
\centerline{\psfig{file=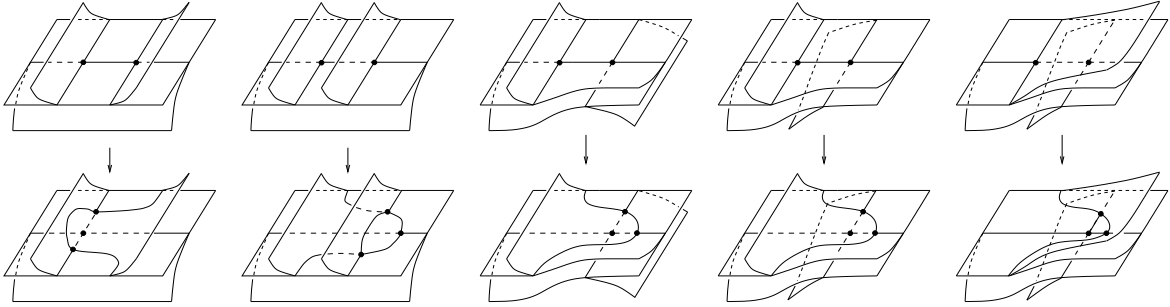,width=15.5truecm}}
\caption{\label{mosse:ramificate}Sliding moves.}
\end{figure}

In~\cite{lnm} we have shown that if $P$ is a branched standard spine and $P'$
is a spine obtained from $P$ by a positive {\rm MP}-move then also $P'$ can be
given the structure of a branched spine (sometimes not in a unique way). Some
of the branched {\rm MP}-moves induce, at least locally, an essential modification of
the black-white splitting of the boundary of the corresponding manifold.
All other {\rm MP}-moves, which we have called sliding moves, do not change up to
isomorphism the pair $(M,v)$ associated to the spine. Moreover the moves can be
realized in $(M,v)$ as continuous modifications of one spine into another
one, through spines positively transversal to $v$, with exactly one
non-standard spine along the modification. If one takes into account
orientations there are exactly 16 different sliding moves, but the essential
phenomena are only those described in Fig.~\ref{mosse:ramificate}.

In~\cite{lnm} we have based on these sliding moves (together with an extra
0-to-2 move) a combinatorial presentation of combed closed manifolds. In
Section~\ref{improve:lnm:section} we will substantially improve this result
to include the case with boundary, and we will also show that the extra move can
be disposed of. This improvement will allow us to define torsions in various
topologically relevant situations, thanks to the following result to which the
present section is devoted:

\begin{teo}\label{sliding:invariance}
All torsions are invariant under sliding moves.
\end{teo}

In view of the warning given in Subsection~\ref{comparison:warning:subsection} the
assertion that torsion is invariant under sliding moves must be interpreted with
some care: as we have remarked above, slidings can be physically realized
as continuous modifications inside combed manifolds. Using such a ``physical
model'' for a sliding move $P_0\move{\ } P_1$, we can consider a ``common model''
for the (abstractly homeomorphic) spaces $X(P_0)$ and $X(P_1)$, so it makes sense
to say that torsions are actually the same.

\dim{sliding:invariance} We start with invariance of the basic version
$\tau^\varphi(P,\hbasis)$ of torsion. As mentioned in the introduction,
our proof mimics the more general proof of Turaev~\cite{turaev:Euler}, but it is
self-contained. The reason why our proof is technically easier is that we only
need to deal with some definite local modifications of the cell complex structure.

We first need to recall the subdivision technique. Consider a branched spine
$P$ and the corresponding complex $X(P)$ with triangulation $\ttt(P)$, and let
$\dd$ be a subdivision of $\ttt(P)$. We will mainly be interested in the case
where also $\dd$ is a triangulation (possibly with multiple and
self-adjacencies). A subdivided spider $s_\dd(P)$ can be defined as
$\sum_{p}\beta_p$, where $\{p\}$ is a collection of one interior point for each
simplex of $\dd$, and $\beta_p$ is (the closure of) the orbit of $v(P)$ which
starts at $p$ and reaches $x_0(P)$. The  reader can easily check that the
choice of $\{p\}$ is inessential, so we omit it from the notation. Our
definition of subdivided spider is somewhat easier than Turaev's general one,
thanks to our choice of taking the simplices of $\ttt(P)$ to be unions of
orbits of $v(P)$. Now consider the initial data 
$\varphi,\Lambda,\hbasis,\orient$ which allow to define a torsion 
$\tau_0^\varphi(P,\hbasis,\orient)$. We can define the $\Lambda$-modules ${\rm
C}^{\varphi,\dd}_*(P)$ using cellular chains with respect to the lifting of
$\dd$ to $\tilde X(P)$, and we can use the spider $s_\dd(P)$ to construct preferred $\Lambda$-bases
of these modules. Now one easily checks that there exists a canonical
isomorphism $H^{\varphi,\dd}_*(P)\cong H^\varphi_*(P)$, so we can use the same
symbol $\hbasis$ for a $\Lambda$-basis of $H^{\varphi,\dd}_*(P)$, and define a
torsion $\tau^{\varphi,\dd}(P,\hbasis)$ exactly as we have done in the previous
section.  The next result is a simplified version of Lemma~3.2.3
in~\cite{turaev:Euler}. 

\begin{prop}\label{subdivision:same}
$\tau^{\varphi,\dd}(P,\hbasis)=\tau^\varphi(P,\hbasis)$.
\end{prop}

\begin{figure}
\centerline{\psfig{file=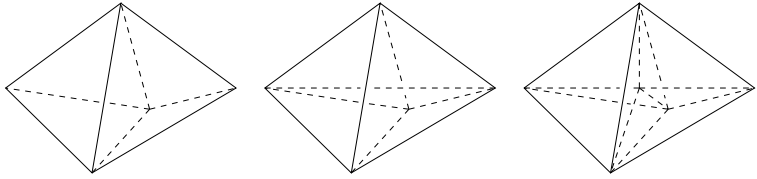,width=10truecm}}
\caption{\label{subdivide} Common subdivision of triangulations.}
\end{figure}

We will prove this proposition only in the more specific
situation we are actually interested in. So, let us consider a sliding move
$P_0\move{\ } P_1$, physically realized inside a manifold $M$, so that the pointed
complexes $X(P_0)$ and $X(P_1)$ are both identified to $X=M/\partial M$. Note that
$X$ comes with two distinct triangulations $\ttt(P_0)$ and $\ttt(P_1)$.
Figure~\ref{subdivide} shows the obvious simplest common subdivision $\dd$ of
$\ttt(P_0)$ and $\ttt(P_1)$. We confine ourselves to showing that 
$\tau^{\varphi,\dd}(P_j,\hbasis)=\tau^\varphi(P_j,\hbasis)$ for $j=0,1$.

\dim{subdivision:same}  Note first that the definition of subdivided spider
leads to the following very natural rule: a simplex of $\dd$ lying in a simplex
$S$ of $\ttt(P)$ is lifted in $\tilde X(P)$ to the only preimage which lies in
the lifting of $S$ determined by $P$. Now, this rule makes sense also for more
general subdivisions than triangulations, in particular for cell complexes, so
we will use them. One easily sees that the subdivisions of $\ttt(P_0)$ and
$\ttt(P_1)$ into $\dd$ can be expressed as combinations of the following
elementary transformations:
\begin{enumerate}
\item An edge is subdivided into two edges by adding a vertex;
\item A square is subdivided into two triangles by adding a diagonal;
\item The inverse of the transformation which removes a triangle, thus replacing
the two adjacent tetrahedra by one polyhedron with 5 vertices, 9 edges, and 6 
triangular faces.
\end{enumerate}
We are left to prove that torsion is invariant under these
transformations. In all three cases, the proof goes as follows: 
\begin{enumerate}
\item[(i)] consider data $\basis_i,\tilde\hbasis_i,g_i$,  $i=0,\dots,3$, which allow to compute $\tau$
before subdivision;
\item[(ii)] describe new data $\basis'_i,\tilde\hbasis'_i,g'_i$
for the subdivided complex;
\item[(iii)] analyze the matrices
$((\partial\basis_{i+1})\cdot\tilde\hbasis_i\cdot\basis_i)/g_i$ and 
$((\partial\basis'_{i+1})\cdot\tilde\hbasis'_i\cdot\basis'_i)/g'_i$
to show that they are the same in $\Kbar_1(\Lambda)$.
\end{enumerate}
Note that this proves that torsion is unchanged ``term by term'', not only
globally.

We only deal with transformation 3, leaving the other cases to the reader. Denote
by $Q$ the polyhedron which is split into tetrahedra $T_1,T_2$ by addition of a
triangle $\Delta$. Then, in a natural way, $g'_0=g_0$, $g'_1=g_1$, and $g'_2$ is
obtained from $g_2$ by adding the lifting of $\Delta$ which lies in the lifting of
$Q$; to get $g'_3$ from $g_3$, we need to remove the lifting of $Q$ and add the
liftings of $T_1$ and $T_2$ which lie in it. For the lifted homology bases, we have
$\tilde\hbasis'_i=\tilde\hbasis_i$ for $i=0,1,2$, and  $\tilde\hbasis'_3$ is
obtained from $\tilde\hbasis_3$ by replacing each occurrence of $\tilde Q$ with
$\tilde T_1+\tilde T_2$. Similarly, $\basis'_i=\basis_i$ except for $i=3$, and
$\basis'_3$ is obtained from $\basis_3$ by replacing each occurrence of $\tilde Q$
with $\tilde T_1+\tilde T_2$, and then adding $\tilde T_2$. The transition matrices
are unchanged in dimensions $0$ and $1$, while in dimensions $2$ and $3$, with
obvious meaning of symbols, we have:
\begin{eqnarray*}
((\partial\basis'_3)\tilde\hbasis'_2\basis'_2)/g'_2
& = &
\pmatrix{ 	& *	&	&	\cr
	\partial\basis_3/g_2 & \vdots & \hbasis_2/g_2 & \basis_2/g_2 \cr
		& *	&	&	\cr
	0\cdots0& 1 	& 0\cdots0&0\cdots0\cr},\\
(\tilde\hbasis'_3\basis'_3)/g'_3
& = &
\pmatrix{ & & 0	\cr
\tilde\hbasis_3/(g_3\setminus\{\tilde Q\}) & \basis_3/(g_3\setminus\{\tilde Q\}) &\vdots \cr
	& 0	\cr
\tilde\hbasis_3/\tilde Q & \basis_3/\tilde Q & 0 \cr
\tilde\hbasis_3/\tilde Q & \basis_3/\tilde Q & 1 \cr}.
\end{eqnarray*}
When $\Lambda$ is a field, one immediately gets the conclusion by
taking determinants. For the general case one needs to recall the definition
of $\Kbar_1(\Lambda)$, and the conclusion follows anyway.
\finedim{subdivision:same}

Our next step is again a simplified version of a more general result. With the
same notation as above, consider the subdivided spiders
$s_\dd(P_i)=\sum_{p}\beta^{(i)}_p$. Using Remark~\ref{non:zero:chi:good:1} it is
easy to see that $\sum_p\varepsilon(p)(\beta^{(0)}_p-\beta^{(1)}_p)$ is a cycle,
where $\varepsilon(p)=(-1)^d$, and $d$ is the dimension of the cell of which $p$
is the centre. We denote by  $h(P_0,P_1)$ the class in $H_1(X(P);\mz)$ of this
cycle.

\begin{prop}\label{compare:slided}
If $h(P_0,P_1)=0$ then $\tau^\varphi(P_0,\hbasis)=
\tau^\varphi(P_1,\hbasis)$.
\end{prop}

\dim{compare:slided} By the previous result, it is enough to show
that $\tau^{\varphi,\dd}(P_0,\hbasis)=\tau^{\varphi,\dd}(P_1,\hbasis)$.
With obvious meaning of symbols, we have
\begin{eqnarray*}
\tau^{\varphi,\dd}(P_1,\hbasis) & = &
\pm\prod_{i=0}^3\Big[\Big((\partial\basis_{i+1})\tilde\hbasis_i\basis_i\Big)
\Big/ g_i^{\varphi,\dd}(P_1)\Big]^{(-1)^i} \\
& = & \pm\prod_{i=0}^3\Big[\Big((\partial\basis_{i+1})\tilde\hbasis_i\basis_i\Big)
\Big/ g_i^{\varphi,\dd}(P_0)\Big]^{(-1)^i}\cdot
\Big[ g_i^{\varphi,\dd}(P_0)
\Big/ g_i^{\varphi,\dd}(P_1)\Big]^{(-1)^i} \\
& = & \tau^{\varphi,\dd}(P_0,\hbasis)\cdot
\Big( \prod_{i=0}^3 \Big[ g_i^{\varphi,\dd}(P_0)
\Big/ g_i^{\varphi,\dd}(P_1)\Big]^{(-1)^i} \Big)
\end{eqnarray*}

To compute the last correction factor, we start by remarking that  the
homomorphism $\varphi:\pi\to\Lambda_*$ induces another one
$\varphi':\pi\to\Kbar_1(\Lambda)$. Since $\Kbar_1(\Lambda)$ is an Abelian group, 
$\varphi'$ induces a homomorphism $\varphi'':H_1(X)\to\Kbar_1(\Lambda)$.
Now, let us denote by $\tilde e_p^{(k)}$ the lifting in $\tilde X$ of the cell of
$\dd$ centred in $p$, where $k=0,1$ and the lifting is determined by $s_\dd(P_k)$.
Note that $\tilde e_p^{(0)}=\gamma_p\cdot\tilde e_p^{(1)}$, with 
$\gamma_p=[(\beta^{(1)}_p)^{-1}\cdot\beta^{(0)}_p]\in\pi_1(X,x_0)$.
This implies that $[g_i^{\varphi,\dd}(P_0)/g_i^{\varphi,\dd}(P_1)]$ is the image in
$\Kbar_1(\Lambda)$ of a diagonal matrix with entries $\varphi(\gamma_p)$, as
$p$ varies in the centres of $i$-cells of $\dd$. It easily follows that
$$\prod_{i=0}^3 \Big[ g_i^{\varphi,\dd}(P_0)
\Big/ g_i^{\varphi,\dd}(P_1)\Big]^{(-1)^i}=\varphi''(h(P_0,P_1))$$
whence the conclusion.
\finedim{compare:slided}
  
The next result concludes the proof of invariance.

\begin{prop}\label{slided:same}
For all sliding moves $P_0\move{\ } P_1$ we have $h(P_0,P_1)=0$.
\end{prop}

\dim{slided:same} Instead of treating in detail all the moves, we confine
ourselves to the description of a general framework, and then we apply the method
to one instance of the move, the other cases being similar. Recall that
Fig.~\ref{subdivide} describes the portion $R$ where the move takes place, with
portions of the triangulations $\ttt(P_0)$, $\ttt(P_1)$ and their common
subdivision $\dd$. To be precise, the figure shows an ``unfolded version'' $R'$ of
the portion $R$, because in $R$ all the ``external'' vertices of the figure are
identified together (with the point $x_0$), so for example each edge of
$\ttt(P_0)$ and $\ttt(P_1)$ represents a (probably non-trivial) element of
$\pi_1(X,x_0)$. Note that also some external edges or faces could be glued together.
The basic idea of the proof is to lift the cycle 
$\sum_p\varepsilon(p)(\beta^{(0)}_p-\beta^{(1)}_p)$ to a 1-chain in $R'$, and show
that the result is again a cycle. Since $R'$ is contractible, the conclusion
easily follows.

\begin{figure} 
\centerline{\psfig{file=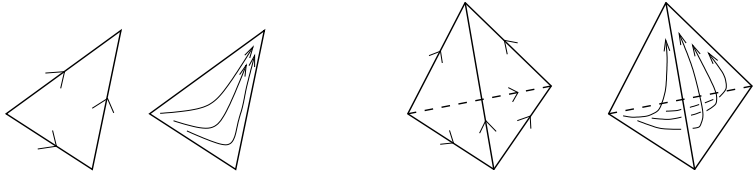,width=10truecm}}
\caption{\label{fieldorbi} How to deduce the field from the edge-orientation.}
\end{figure}

We remark now that a branching on a standard spine can be encoded just by an
orientation of the edges of the dual triangulation. This orientation is defined
by the requirement that the algebraic intersection between a region of the
spine and the dual edge should be positive. Edge-orientations coming from
branchings are characterized by the property that no triangle has a circular
orientation of its edges (see also~\cite{gr:1}). Recall as well that the edges
are orbits of the field defined by the spine, and orientations match. More
importantly, edge-orientations  allow us to describe the orbits of the field
also on the open triangles and tetrahedra of the triangulation. The rules are
illustrated in Fig.~\ref{fieldorbi}. 

Now, each sliding move will be encoded by a pair of matching patterns of
orientations of the edges of $\ttt(P_0)$ and $\ttt(P_1)$, both defining fields
on the common portion they triangulate. Note that non-zero contributions to
$h(P_0,P_1)$ can only come from simplices of $\dd$ which are not shared with
both $\ttt(P_0)$ and $\ttt(P_1)$. This rules out all simplices outside the
portion $R'$ and those on its boundary. We are left to deal with internal
simplices, namely 1 vertex, 5 edges, 9 faces and 6 tetrahedra. To carry out the
program outlined above, we must lift to $R'$ the corresponding loops
$(-1)^d\cdot(\beta^{(0)}_p-\beta^{(1)}_p)$ and show that the sum of the
boundaries is always null.

\begin{figure} 
\centerline{\psfig{file=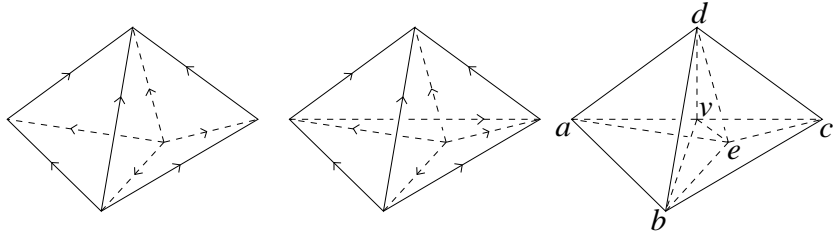,width=11truecm}}
\caption{\label{calcex1} A sliding move, and notations for the subdivided
triangulation.}
\end{figure}

As promised, we carry out calculations in one example. The first sliding move
of Fig.~\ref{mosse:ramificate} translates in terms of edge-orientations
as described in Fig.~\ref{calcex1},
where again we are showing the ``unfolded version''.
In the same figure we introduce some useful notation for the subdivided
triangulation. Note again that $a=b=c=d=e=x_0$ in $X$.
To analyze the contributions of the internal simplices to the boundary
of the lifted chain, we need to determine, for both fields, the
targets of the orbits starting at centres of simplices. This is done in the next
tables.
\vspace{.2cm}

{\small

\begin{center}

\begin{tabular}{||l||c|c|c|c|c|c|c|c|c|c|}
\hline\hline
Simplex $\sigma$ & 
$v$ & $va$ & $vb$ & $vc$ &  $vd$ & $ve$ &
$vab$ & $vad$ & $vae$ & $vcb$ 
\\ \hline
$(-1)^{{\rm dim}(\sigma)}$ &
+1 & $-1$ & $-1$ & $-1$ & $-1$ & $-1$ & +1 & +1 & +1 & +1  
\\ \hline
End $\tilde\beta^{(0)}_{p(\sigma)}$ & 
$d$ & $d$ & $d$ & $d$ & $d$ & $d$ & $d$ & $d$ & $d$ & $d$  
\\ \hline
End $\tilde\beta^{(1)}_{p(\sigma)}$ & 
$c$ & $c$ & $c$ & $c$ & $d$ & $c$ & $c$ & $d$ & $c$ & $c$  
\\ \hline
Boundary & 
$d-c$ & $c-d$ &  $c-d$ &  $c-d$ & $0$ & $c-d$ &
$d-c$ & 0 & $d-c$ & $d-c$ 
\\ \hline\hline
\end{tabular}

\vspace{.3cm}

\begin{tabular}{||l||c|c|c|c|c|c|c|c|c|c|c|}
\hline\hline
Simplex $\sigma$ & 
$vcd$ & $vce$ & $vbe$ & $ved$ & $vdb$ & $vabe$ & $vaed$ & $vadb$ & $vcbe$ & $vced$ & $vcdb$
\\ \hline
$(-1)^{{\rm dim}(\sigma)}$ &
+1 & +1 & +1 & +1 & +1 & $-1$ & $-1$ & $-1$ & $-1$ & $-1$ & $-1$ 
\\ \hline
End $\tilde\beta^{(0)}_{p(\sigma)}$ & 
$d$ & $d$ & $d$ & $d$ & $d$ & $d$ & $d$ & $d$ & $d$ & $d$ & $d$
\\ \hline
End $\tilde\beta^{(1)}_{p(\sigma)}$ & 
$d$ & $c$ & $c$ & $d$ & $d$ & $c$ & $d$ & $d$ & $c$ & $d$ & $d$ 
\\ \hline
Boundary & 
$0$ & $d-c$ & $d-c$ & $0$ & $0$ & $c-d$ & $0$ & $0$ & $c-d$ & $0$ & $0$ 
\\ \hline\hline
\end{tabular}

\end{center}

}

\vspace{.2cm}

\noindent The sum of the bottom rows of the two tables is null, which gives the
desired conclusion. This eventually proves the invariance of 
$\tau^\varphi(P,\hbasis)$. The sign-refined version
$\tau_0^\varphi(P,\hbasis,\orient)$ is dealt with analogously, with a little more
effort. A similar argument works  also for  the modified versions of torsion
$\tau_{\rm w}$ and $\tau_{{\rm w},0}$. One only needs to note that not all the vertices which appear
Fig.~\ref{subdivide} are identified to the basepoint $x_0$, but the vertices which
are endpoints of orbits indeed are, so the proof proceeds
exactly as above.\finedim{sliding:invariance}

\subsection{The canonical spider and the Euler class}

Let us recall that our definition of torsion of a spine was based on a spider
obtained by integrating the field in the {\it positive} direction starting from
the centres of the cells. However, as the reader can easily check, another
definition is obtained along the same lines but integrating the field in the {\em
negative} rather than positive direction. The problem naturally arises to compare
these constructions. We will now prove that the difference of the positive and the
negative spiders, when lifted to $M(P)$, is a very natural object, namely the
Poincar\'e dual of the Euler class of the plane distribution normal to the field.
See Subsection~\ref{discussion} for further comments on the meaning of this result.

For a precise statement, we need to introduce some notation. Given a branched
spine $P$, we denote by $\sum_c\tilde\alpha^{+}_c$ the natural lifting to
$M(P)$ of the spider defined in Subsection~\ref{basic:subsection}. Note that
this is not quite a spider any more, because the head is exploded into a set of
points scattered on $\partial(M(P))$. We define $\sum_c\tilde\alpha^{-}_c$ in a
similar way, integrating the field in the negative direction. So, each
difference $\tilde\alpha^{+}_c-\tilde\alpha^{-}_c$ represents an arc with both
ends on $\partial(M(P))$. Now, with the usual meaning of symbols, we have:

\begin{prop}\label{Euler:class:prop}
$[\sum_c\varepsilon(c)\cdot(\tilde\alpha^{+}_c-\tilde\alpha^{-})]=
{\rm PD}(\ee(v(P)^\perp))\in H_1(M(P),\partial(M(P));\mz)$.
\end{prop}

\begin{figure}
\centerline{\psfig{file=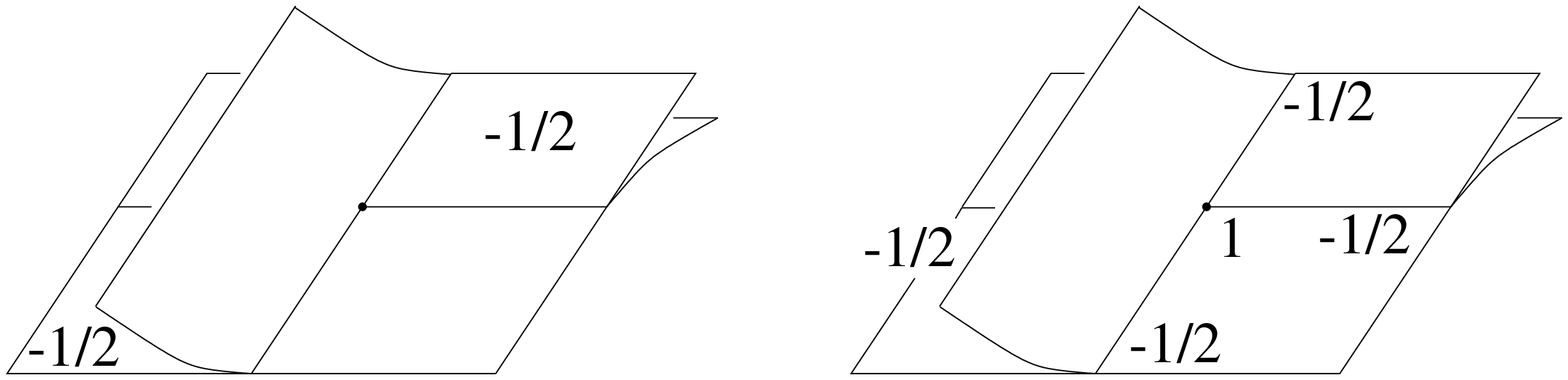,width=9cm}}
\caption{\label{maw:contributions} Contributions of a vertex to the cochains.}
\end{figure}

\dim{Euler:class:prop} In~\cite{lnm}, inspired by~\cite{christy}, we have
introduced a canonical cochain $c_P$ for the Euler class of a field encoded by a
branched spine $P$. To get $c_P$, one considers on $P$, near the singular set
$S(P)$, the tangent field $\mu_P$ which is transversal to $S(P)$ and points
from the locally two-sheeted area to the locally one-sheeted area. The
value of $c_P$ on a region $R$ is the index of the extension of $\mu_P$ to $R$.
Noting that, at each vertex, $\mu_P$ is tangent to the boundary of exactly two
(opposite) regions, one sees that $c_P(R)=1-n(R)/2$, where $n(R)$ is the number
(with multiplicity) of vertices of $R$ at which $\mu_P$ is tangent to $\partial
R$. Contributions of a vertex to the regions incident to it are shown on the
left in Fig.~\ref{maw:contributions}.

Now, in $\sum_c\varepsilon(c)\cdot(\tilde\alpha^{+}_c-\tilde\alpha^{-})$, 
arcs coming from vertices and edges of $P$ are not in
general position with respect to $P$, so we modify them by slightly pushing
along $\mu_P$. The value on a region $R$ is now given by a sum of
contributions: one +1, coming from $R$ itself, some $-1$'s, coming from edges,
and some $+1$'s, coming from vertices. If one halves the contributions of edges
and localizes the halves at the ends, one has that the values of the cochain
are expressed just as for $c_P$, namely $+1$ plus some contributions coming
from vertices. The latter are shown on the right in
Fig.~\ref{maw:contributions}. Since contributions are actually the same, the
proof is complete.\finedim{Euler:class:prop}

\section{Branched spines of combed manifolds\\ with 
boundary pattern}\label{improve:lnm:section}

In this section we will extend the main results of Chapter~5 of~\cite{lnm} from
the closed to the bounded case. Namely, we will show that compact manifolds
with concave combings (see below for the precise definitions) are
combinatorially described by (suitable) branched spines up to sliding. We first
show that one of the sliding moves considered in~\cite{lnm} is essentially
generated by the other moves. 

\subsection{The snake move} 
\begin{figure}
\centerline{\psfig{file=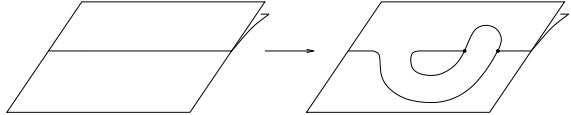,width=7.5cm}}
\caption{\label{snake:move} The snake move.}
\end{figure}

In
Section~\ref{sliding:invariance:section}, when proving invariance of torsions,
we have not dealt with the extra
move which, together with the branched {\rm MP}-moves, was defined in~\cite{lnm} to be a
generator of the sliding calculus. This is the ``snake move'', described in
Fig.~\ref{snake:move}
(actually, taking into account orientations, there are two versions of the
move). Our reason for not treating the move was the next
Proposition~\ref{rigid:spines} (see also 
Corollary~\ref{no:snake:nec}). 

In the rest of this section we will denote by $\pp$ the set of (isomorphism
classes of) branched spines, and by $\rr$ the subset of those
which are ``rigid'' from the point of view of the branched {\rm MP}-moves, {\it
i.e.}~the spines to which no such move applies. An explicit description of $\rr$
is given in the proof of the next result. In the statement we only emphasize the
most important consequences of this description. It will be convenient to use the
terminology introduced in~\cite{lnm}: we call {\em trivial} and denote by 
$\stwotriv$ the bicoloration of $S^2$ which consists of one black disc and one
white disc.

\begin{prop}\label{rigid:spines}
\begin{enumerate}
\item[(i)] For every bicolorated surface $\Sigma$ there are at most two spines
$P\in\rr$ such that $\partial(M(P))\cong\Sigma$.
\item[(ii)] If two elements of $\pp\setminus\rr$ are related through
branched {\rm MP}-moves and snake moves, they are also related through 
branched {\rm MP}-moves only.
\end{enumerate}
\end{prop}

\begin{figure}
\centerline{\psfig{file=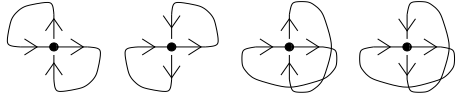,width=6cm}}
\caption{\label{one:vert:spi} Branched spines with one vertex.}
\end{figure}

\dim{rigid:spines} In this proof we will find convenient to use the graphic
representation of branched spines introduced in~\cite{lnm}. We will skip most
details, giving only the main points. 
We start by listing rigid spines. Note first that if a negative branched
{\rm MP}-move applies to a spine then also a positive one does, so we only need
to consider positive rigidity. The spines with one vertex, shown in
Fig.~\ref{one:vert:spi}, are of course rigid.
Using~\cite{lnm} one easily checks that $\partial(M(P))$ is $\stwotriv$ for the
first two spines, and $\stwotriv\sqcup\stwotriv$ for the other two.

\begin{figure}
\centerline{\psfig{file=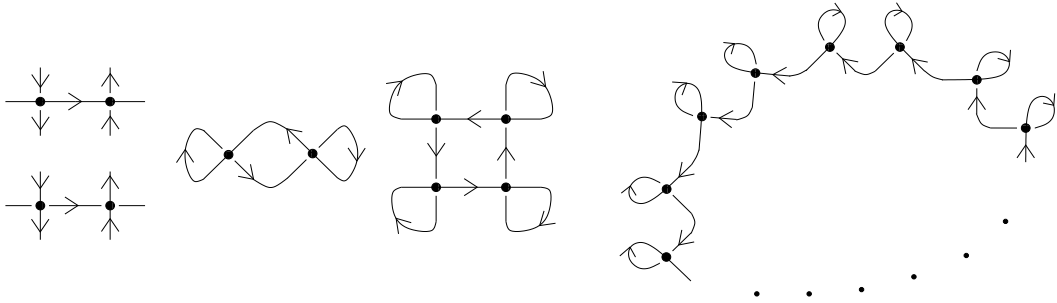,width=14cm}}
\caption{\label{rigid:spines:fig} Rigid branched spines.}
\end{figure}

Now we turn to rigid spines with more than one vertex. Rigidity implies that
all edges with distinct endpoints should appear as on the left in 
Fig.~\ref{rigid:spines:fig}.
It is not hard to deduce that rigid spines come in a sequence $P_1^{\rm rig}$,
$P_2^{\rm rig}$, $\dots$ as
shown in the rest of Fig.~\ref{rigid:spines:fig}, where $P_k^{\rm rig}$
has $2k$ vertices, and $\partial(M(P_k^{\rm rig}))$ is the union of
$\stwotriv$ together with $k$ copies of $S^2_{\rm white}$ and 
$k$ copies of $S^2_{\rm black}$.
This classification proves {\em (i)}. 

To show {\em (ii)} we must prove that: 
\begin{enumerate}
\item[{\em (ii-a)}] Sequences which contain rigid spines can be replaced
by sequences which do not.
\item[{\em (ii-b)}] If two non-rigid spines are related by
one snake move then they are also related by a sequence of branched {\rm MP}-moves.
\end{enumerate}

For {\em (ii-a)}, we note that the result of a positive snake move is never
rigid. So if a rigid spine $P$ appears in a sequence of moves then around $P$
we see $P_1\move{\mu_1^{-1}}P\move{\mu_2}P_2$, with $\mu_1$ and $\mu_2$
positive snake moves. Since all edges of a spine survive through a snake move,
there is a version $\tilde\mu_2$ of $\mu_2$ which applies to $P_1$ and 
a version $\tilde\mu_1$ of $\mu_1$ which
applies to $P_2$, and the result is the same. So we have
$P_1\move{\tilde\mu_2}\tilde P\move{\tilde\mu_1^{-1}}P_2$, and now all the
spines involved are non-rigid.  

\begin{figure}
\centerline{\psfig{file=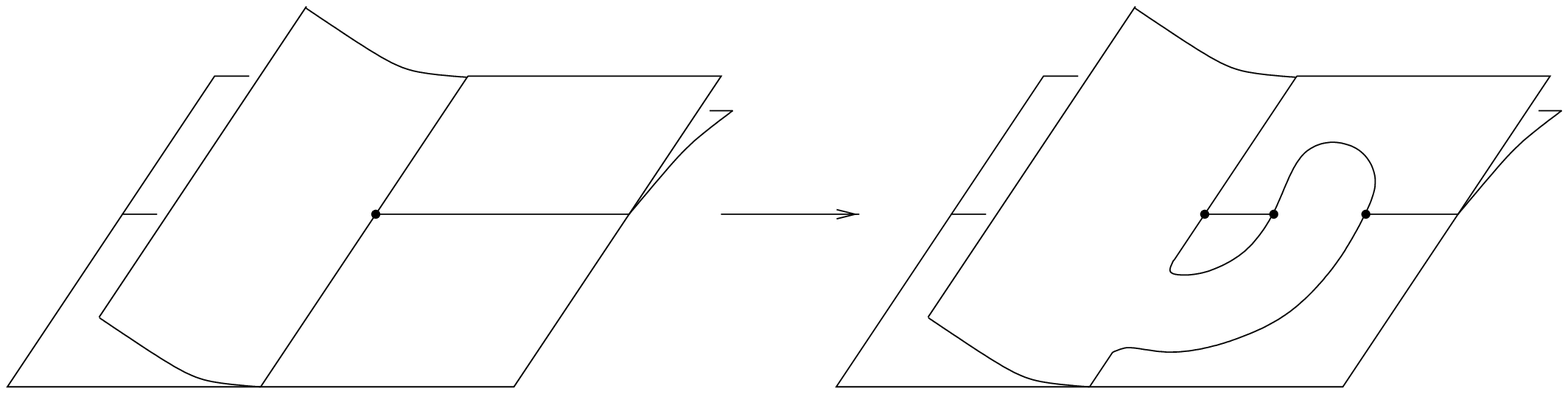,width=7.5cm}}
\caption{\label{vertex:move} The vertex move.}
\end{figure}

Let us turn to {\em (ii-b)}. The proof results from three steps, to describe which
we introduce in Figure~\ref{vertex:move}  
another move, called vertex move, whose unbranched version was already
considered in~\cite{matv:mossa} and~\cite{piergallini}. Again, taking into
account orientations, there are two versions of the move (for each vertex
type), but we will ignore this detail. 

{\em Step 1: if $v$ is a vertex of a spine $P$, $e$ is any one of the
edges incident to $v$, $P_v$ is obtained from $P$ via the vertex move at $v$,
and $P_e$ is obtained from $P$ via the snake move on $e$, then $P_v$ and
$P_e$ are related by {\rm MP}-moves.} This is proved by an easy case-by-case
analysis. It turns out that two {\rm MP}-moves (a positive and a negative one)
are always sufficient.

{\em Step 2: let $v$, $P$ and $P_v$ be as above. If $P$ and $P_v$  are related
by {\rm MP}-moves, the same is true for $P$ and any spine obtained from $P$ by
a snake move.} To see this, use step 1 to successively transform vertex moves
into snake moves and conversely, until the desired snake move is reached.

{\em Step 3: if $P$ is non-rigid then there exists a vertex $v$ such that
$P$ and $P_v$ are related by {\rm MP}-moves.} The vertex $v$ is chosen to be an
endpoint of an edge to which the positive {\rm MP}-move applies. The argument is
again a long case-by-case one, which refines in a branched context the argument
given by Piergallini in~\cite{piergallini}. The sequence always consists of
three positive moves followed by a negative one.
This concludes the proof.\finedim{rigid:spines}

\subsection{Generalized combed calculus}

Recall from Subsection~\ref{reminder} that a {\em concave} field $v$ on a
compact 3-manifold $M$ is one which is tangent to $\partial M$ only ``from
inside'', along some simple curves. We will denote by $\combd$ the set of all
such pairs $(M,v)$, viewed up to homeomorphism of $M$ and homotopy of $v$
through concave fields. Note that the black-white splitting of $\partial M$
evolves isotopically during a homotopy of $v$, so we can associate to
$[M,v]\in\combd$ a well-defined {\em boundary pattern}. A class
$[M,v]\in\combd$ is called a {\it combing} on the homeomorphism class of the
manifold $M$. For a technical reason we rule out from $\combd$ the set of those
classes $[M,v]$ such that $\partial M$ contains components of the type
$\stwotriv$. This is actually not a serious restriction, because each
$\stwotriv$ component can be capped off by a $\bthtr$ (the 3-ball with constant
field), and the result is well-defined  up to homotopy. Note that we do accept
pairs $(M,v)$ with $M$ closed, and pairs in which $v$ has no tangency at all to
$\partial M$.

Let us denote now by $\vv$ the set of pairs $(M,v)$ where $v$ is concave and
{\em traversing}, {\it i.e.}~such that all orbits are segments with both ends
on $\partial M$, and in $\partial M$ there is exactly one $\stwotriv$
component. These pairs will be viewed up to homeomorphism of $M$ and
homotopies of $v$ through concave traversing fields. By the
construction recalled in Subsection~\ref{reminder}, we can associate to every
branched spine a manifold with a concave traversing field. We will denote by
$\objd$ the set of (isomorphism classes of) those spines which give rise to
elements of $\vv$. Given $[P]\in\objd$, if we cap off the only $\stwotriv$ in
$\partial(M(P))$ by a $\bthtr$, we obtain a well-defined element of $\combd$.

\begin{teo}\label{main:improvement:teo}
The map $\recd:\objd\to\combd$ thus defined is surjective, and the equivalence
relation defined by $\recd$ on $\objd$ is generated by standard sliding moves.
\end{teo}

\begin{rem}\label{trivial:pieces}
{\em The following interpretation of the surjectivity of $\recd$ is perhaps useful.
Note first that the dynamics of a field, even a concave one, can be very complicated,
whereas the dynamics of a traversing field (in particular, $B^3_{\rm triv}$) is simple.
Surjectivity of $\recd$ means that for any (complicated) concave field
there exists a sphere $S^2$ which splits the field into two (simple) pieces:
a standard $B^3_{\rm triv}$ and a concave traversing field. Actually, a 1-parameter version
of this statement also holds (see Remark~\ref{one-parameter}): we will need it to show that the 
$\recd$-equivalence is the same as the sliding equivalence.}
\end{rem}

\begin{rem}\label{no:snake:nec}
{\em By Proposition~\ref{rigid:spines}, in Theorem~\ref{main:improvement:teo}
we could remove from $\objd$ the set $\rr$
of {\rm MP}-rigid spines and forget the snake move, 
leaving the rest of the statement unchanged.}
\end{rem}

The proof of Theorem~\ref{main:improvement:teo} is an extension of the argument
given in Chapter~5 of~\cite{lnm}, and it is based on the following technical
notion, which extends ideas originally due to Ishii~\cite{ishii}. Let $v$ be a
concave field on $M$. Let $B_1,\dots,B_k$ be the black components of the
splitting of $\partial M$, {\it i.e.}~the regions on which $v$ points outwards. A
{\em normal section} for $(M,v)$ is a compact surface $\Sigma$ with boundary,
embedded in the interior of $M$, with the following properties:
\begin{enumerate}
\item\label{transverse:point} $v$ is transverse to $\Sigma$;
\item\label{disc:point} $\Sigma$ has exactly $k+1$ components $\Sigma_0,\dots,\Sigma_k$, with 
$\Sigma_0\cong D^2$;
\item\label{black:point} For $i>0$, the projection of $B_i$ on $\Sigma$ along the orbits of $-v$ is
well-defined and yields a homeomorphism between $B_i$ and a surface $B'_i$
contained in the interior of $\Sigma_i$, with $\Sigma_i\setminus B'_i$ being a
collar on $\partial\Sigma_i$;
\item\label{all:orbits:met} Each positive half-orbit of $v$ meets either the interior of some $B_i$
(where it stops), or the interior of some $\Sigma_i$;
\item\label{generic:point} $\partial\Sigma$ meets itself generically along $v$ 
({\it i.e.}~each orbit of
$v$ meets $\Sigma$ at most two consecutive times on $\partial\Sigma$, and, if
so, transversely);
\item\label{standardness:point} Let $P_\Sigma$ be the union of $\Sigma$ with all the orbit
segments starting on $\partial\Sigma$ and ending on $\Sigma$. Then $\Sigma$,
which is a quasi-standard polyhedron by the previous point, is actually
standard.
\end{enumerate}

The next two lemmas show that normal sections of $(M,v)$ correspond bijectively
to spines $P$ such that $\recd(P)=[M,v]$. The proof of surjectivity of $\recd$
and the discussion of its non-injectivity will be based on these lemmas.

\begin{lem}\label{section:then:spine}
If $(M,v)$, $\Sigma$ and $P_\Sigma$ are as above, then $P_\Sigma$ can be given
a structure of branched spine such that $\recd([P_\Sigma])=[M,v]$.
\end{lem}

\dim{section:then:spine}
We orient $\Sigma$ so that $v\trasvint\!^+\,\Sigma$
(by default $M$ is oriented). Every region of $P_\Sigma$
contains some open portion of $\Sigma$, so it can be oriented canonically.
With the obvious screw-orientation, this turns $P_\Sigma$ into a branched spine.

We show that $\recd([P_\Sigma])=[M,v]$ by embedding the abstract  
manifold $M(P_\Sigma)$ in $M$, in such a way that the field carried by
$P_\Sigma$ on $M(P_\Sigma)\subset M$ is just the restriction of $v$. By
construction, $M\setminus M(P_\Sigma)$ will
consist of a copy of $\bthtr$, together with a collar on $\partial M$ which can
be parametrized as $(\partial M)\times[0,1]$ in such a way that $v$ is constant
in the $[0,1]$-direction. This easily implies that $\recd([P_\Sigma])=[M,v]$
indeed.

\begin{figure}
\centerline{\psfig{file=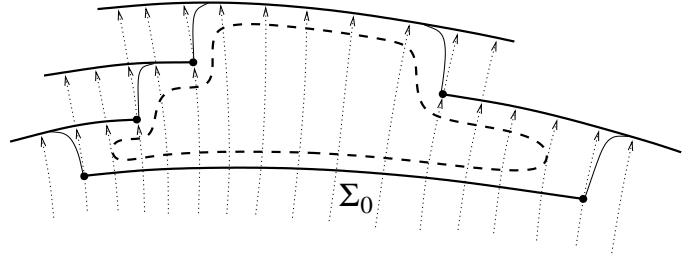,width=9cm}}
\caption{\label{trivball} The trivial ball.}
\end{figure}

We illustrate the embedding of $M(P_\Sigma)$ in $M$ 
pictorially in one dimension less.
Figure~\ref{trivball} shows how $\Sigma_0$ gives rise to a $\bthtr$. 
In the figure we describe $v$ by dotted lines,
$\Sigma$ by thick lines, portions of $P_\Sigma\setminus\Sigma$ by thin lines,
and $\partial(M(P_\Sigma))$ by a thick dashed line. Note also that the portions
of $P_\Sigma\setminus\Sigma$ have been slightly modified so to become positively
transversal to $v$, which allows us to represent the branching as usual,
{\it i.e.}~as a $\cont^1$ structure on $P_\Sigma$.

\begin{figure}
\centerline{\psfig{file=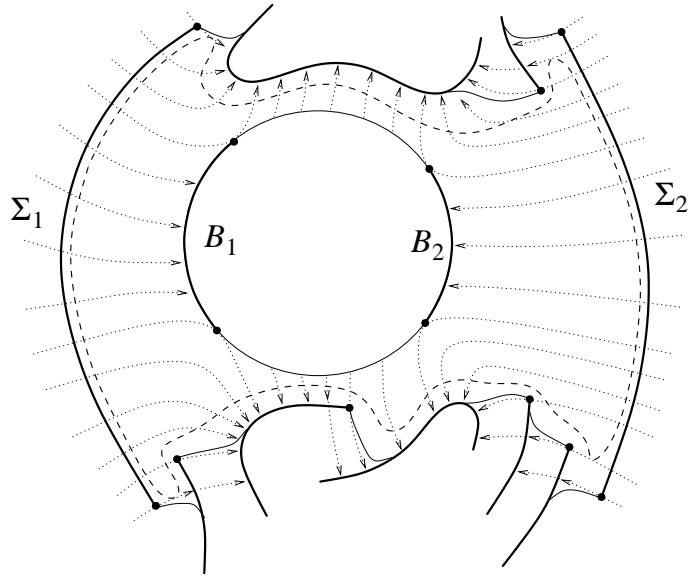,width=9cm}}
\caption{\label{collar} Collar on a boundary component.}
\end{figure}

Figure~\ref{collar} shows the collar based on a component of $\partial M$.
We use the same conventions as in the previous figure, and in addition we
represent the black and white components of $\partial M$ by thick and thin lines
respectively. This description concludes the proof.
\finedim{section:then:spine}

\begin{lem}\label{spine:then:section}
Let $[P]\in\objd$ and $\recd([P])=[M,v]\in\combd$, with $P$ embedded in 
$(M,v)$ according to the geometric description of $\recd$. Let $\Sigma$ be
obtained from $P$ as suggested (in one dimension less) in Fig.~\ref{cutspine}.
\begin{figure}
\centerline{\psfig{file=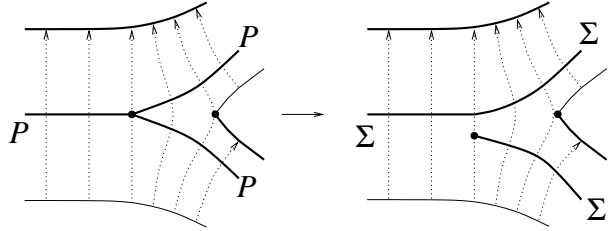,width=8cm}}
\caption{\label{cutspine} Normal section from a spine.}
\end{figure}
Then $\Sigma$ is a normal section of $(M,v)$, and $P_\Sigma$ is isomorphic to
$\Sigma$.
\end{lem}

\dim{spine:then:section}
The construction suggested by Fig.~\ref{cutspine} is obviously the inverse of 
the construction in the proof of Lemma~\ref{section:then:spine}.\finedim{spine:then:section}

\dim{main:improvement:teo} We start with the proof of surjectivity. So, let
us consider a combed manifold $(M,v)$, subject to the usual restrictions. By
Lemma~\ref{section:then:spine} it is natural to try and construct a normal
section for $(M,v)$. Let $B_1,\dots,B_k$ be the black regions in $\partial M$.
Slightly translate each $B_i$ along $-v$, getting $B'_i$. Add to each $B'_i$ a
small collar normal to $v$, getting $\Sigma_i$ (if $\partial B_i=\emptyset$,
we set $\Sigma_i=B'_i$). Select finitely many discs $\{D_n\}$ disjoint from
each other and from all the $\Sigma_i$'s, such that all positive orbits of $v$,
except for the small segments between $B'_i$ and $B_i$, meet
$(\bigcup_{i\geq1}\Sigma_i)\cup(\bigcup D_n)$ in some interior point. Connect
the $D_n$'s together by strips normal to $v$ and disjoint from
$\bigcup_{i\geq1}\Sigma_i$, getting a disc $\Sigma_0$. Up to a generic small
perturbation, the surface $\Sigma=\bigcup_{i\geq0}\Sigma_i$ satisfies all
axioms of a normal section for $(M,v)$, except axiom~\ref{standardness:point}.

Now, even if it is not standard, $P_\Sigma$ can be defined, and the proof of
Lemma~\ref{section:then:spine} shows that it is a quasi-standard branched spine
of $(M\setminus B^3,v)$. In particular, $P_\Sigma$ is connected and its
singular locus is non-empty. Under these assumptions, it is not too hard to see
that there exists a sequence of (abstract) quasi-standard sliding moves which
turns $P_\Sigma$ into a standard spine. If we physically realize these moves
within $M$, preserving transversality to $v$, the result is a standard branched
spine $P$ such that $\recd([P])=[M,v]$.

We are left to show that if $\recd([P_0])=\recd([P_1])$ then $P_0$ and $P_1$
are sliding-equivalent. By the definition of $\combd$ and $\recd$, using also
the above lemmas, there exists a manifold $M$ and a homotopy $(v_t)$ of concave
fields on $M$, such that $P_0$ and $P_1$ are defined by normal sections
$\Sigma^{(0)}$ and $\Sigma^{(1)}$ of $(M,v_0)$ and $(M,v_1)$ respectively.

We prove that $P_0$ and $P_1$ are sliding-equivalent first in the special case
where $v_0=v_1=v$. The general case will be an easy consequence. For $j=0,1$,
let $\Sigma^{(j)}=\bigcup_{i\geq0}\Sigma^{(j)}_i$. Proceeding as in the above
proof of surjectivity, for each black region $B_i$ of $\partial M$, we consider
a collared negative translate $\overline\Sigma_i$ of $B_i$. We choose
$\overline\Sigma_i$ so close to $B_i$ that 
$\overline\Sigma_i\cap\Sigma^{(j)}=\emptyset$, and the negative integration of
$v$ yields a homeomorphism from $\overline\Sigma_i$ to a subset of
$\Sigma^{(j)}_i$. 

\begin{figure}
\centerline{\psfig{file=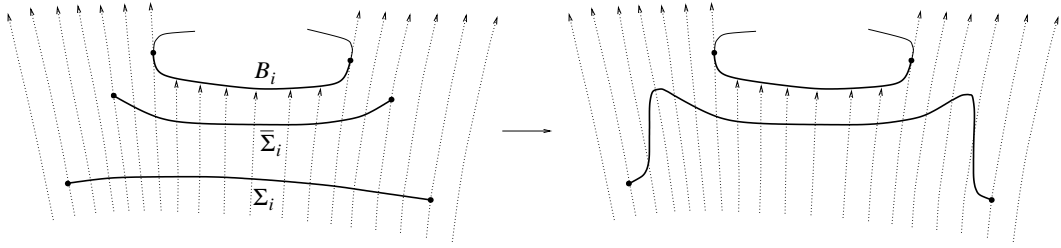,width=14cm}}
\caption{\label{modify:1} Isotopy of a normal section.}
\end{figure}

{\sc Step I}. {\em For $j=0,1$, there exists a disc $D_j$ such that
$D_j\cup(\bigcup_{i\geq1}\overline\Sigma_i)$ is a normal section of $(M,v)$, and
the associated branched spine is sliding-equivalent to $P_j$.} To prove this, we temporarily
drop the index $j$. We first isotope each
$\Sigma_i$, without changing the associated spine, until it contains
$\overline\Sigma_i$, as suggested in Fig.~\ref{modify:1}.

\begin{figure}
\centerline{\psfig{file=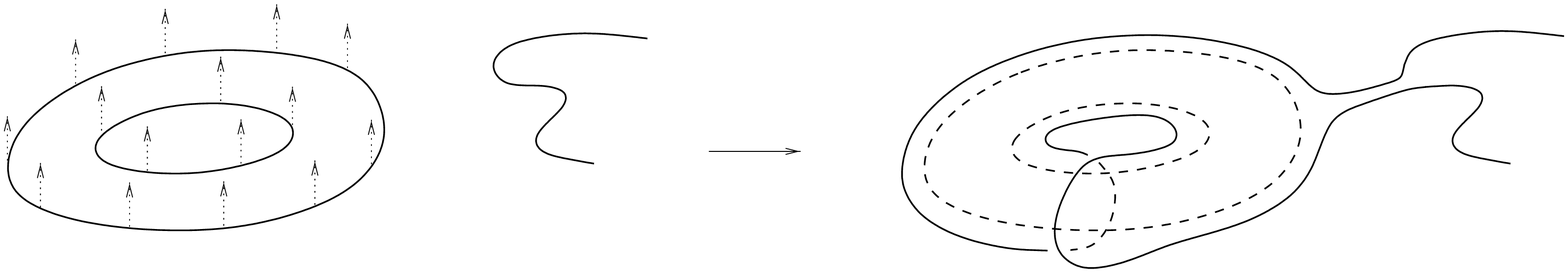,width=14cm}}
\caption{\label{modify:2} Insulation of an annulus.}
\end{figure}

Note that if $\partial B_i=\emptyset$ we automatically have 
$\Sigma_i=\overline\Sigma_i$. Otherwise, we concentrate on one of the
annuli $A$ of which $\Sigma_i\setminus\overline\Sigma_i$ consists. Note
that we cannot just shrink $A$ leaving the rest of the section
unchanged, because we could spoil axiom~\ref{all:orbits:met} of the definition
of normal section. To actually shrink $A$ we first need to ``insulate''
it, toward the positive direction of $v$, by adding to the disc $\Sigma_0$
a strip normal to $v$. Figure~\ref{modify:2} suggests how to do this.

As we modify $\Sigma_0$ as suggested, it is clear that we keep having a
``quasi-normal'' section, {\em i.e.}~all axioms except~\ref{standardness:point}
hold. Moreover the corresponding quasi-standard branched spines are obtained
from each other by quasi-standard sliding moves. To conclude we apply, as
above, the fact that a quasi-standard branched spine is sliding-equivalent to a
standard one, and the technical result established in~\cite{lnm}, according to
which standard spines which are equivalent under quasi-standard sliding moves
are also equivalent under standard sliding moves. This proves Step I.

The conclusion will now follow quite closely the argument in~\cite{lnm}.

{\sc Step II}. {\em There exist discs $D'_0$ and $D'_1$ such that 
$D'_j\cup(\bigcup_{i\geq1}\overline\Sigma_i)$ is a normal section of $(M,v)$
for $j=0,1$, and $D_0\cap D'_0=D'_0\cap D'_1=D'_1\cap D_1=\emptyset$.} Choosing
a metric on $M$, one can construct $D'_0$ and $D'_1$ by first taking many very
small discs almost orthogonal to $v$, and then connecting these discs by strips
transversal to $v$. 

\begin{figure}
\centerline{\psfig{file=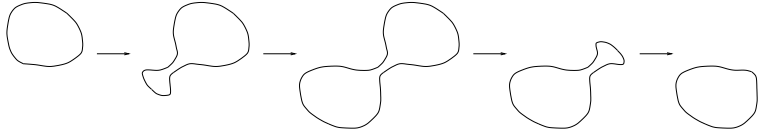,width=10cm}}
\caption{\label{empty:fill} Transformation of a disc into a disjoint one.}
\end{figure}

{\sc Step III}. {\em Conclusion in the case $v_0=v_1$}. If we connect $D_0$ and
$D'_0$ by a strip orthogonal to $v$, we get a bigger disc $\tilde D_0$ such
that $\tilde D_0\cup(\bigcup_{i\geq1}\overline\Sigma_i)$ is still a normal
section of $(M,v)$. We can actually imagine a dynamical process, in which $D_0$
is first enlarged to $\tilde D_0$, and then is reduced to $D'_0$, as in
Fig.~\ref{empty:fill}.
If the transformation is chosen generic enough, at all times
axioms~\ref{transverse:point},~\ref{disc:point},~\ref{black:point} 
and~\ref{all:orbits:met} will hold, and axiom~\ref{generic:point} will
hold at all but finitely many times. This means that the corresponding branched
spines are related by quasi-standard sliding moves. Similarly, we can replace 
$D'_0$ first by $D'_1$ and then by $D_1$.
Using the facts quoted
above, the conclusion follows.

We are left to deal with the general case, where $(v_t)$ is a non-constant
homotopy. It is then sufficient to take a partition $0=t_0<t_1<\dots<t_n=1$ of
$[0,1]$, fine enough that $(M,v_{t_{k-1}})$ and $(M,v_{t_{k}})$ admit a common
normal section which gives rise to isomorphic branched
spines.\finedim{main:improvement:teo}

\begin{rem}\label{one-parameter}
{\em Along the lines of the previous proof we have established
the following topological fact, whose statement does not involve spines. Let
$(v_t)$ be a homotopy of concave fields on $M$, let $B_0,B_1\subset M$ be balls with
$(B_j,v_j)\cong\bthtr$ and $v_j$ traversing on $M\setminus B_j$ for $j=0,1$.
Then there exist another homotopy $(v'_t)$ between $v_0$ and $v_1$ and an
isotopy $(B_t)$ with $(B_t,v_t)\cong\bthtr$ and $v_t$ traversing on $M\setminus
B_t$ for all $t$.}\end{rem}

\section{Torsions of combings and links}

In this section we will combine the results previously obtained, to define
torsions in some topologically relevant situations. We will also give some
hints on how to carry out computations.

\subsection{Torsion of a combed manifold}\label{discussion} As a direct
application of what we have proved, we see that for $[M,v]\in\combd$, if we
choose $[P]\in(\recd)^{-1}([M,v])$, a representation
$\varphi:\pi_1(X(P))\to\Lambda_*$, a $\Lambda$-basis $\hbasis$ of
$H^\varphi_*(P)$ (assuming it to be free) and a homological orientation $\orient$ of $X(P)$, then we
have a torsion $\tau_0^\varphi(P,\hbasis,\orient)\in K_1(\Lambda)$ whose
equivalence class under the natural action of the mapping class group of $M$
depends on $[M,v]$, $\hbasis$ and $\orient$ only.  If the homological orientation is omitted, the
invariant takes values in $\Kbar_1(\Lambda)$. When in $(M,v)$ there are white
boundary components, the definition of $\tau_{{\rm
w},0}^\varphi(P,\hbasis,\orient)$ also applies. 

\begin{rem}\label{one:manifold}{\em The reason why we are forced to let mapping class groups act on
torsions is that in the definition of $\combd$ we have considered manifolds to
be defined only up to homeomorphism. Focusing on a certain marked manifold one
can neglect the action, obtaining a Teichm\"uller-type theory of torsions.
Namely, if a manifold $M$ and a black-white boundary pattern $\pp$ are fixed,
we can consider the set ${\rm Comb}(M,\pp)$ of concave  non-singular vector
fields on $M$ matching $\pp$, viewed up to homotopy. Now $X=M/(\partial
M\cup\{*\})$, and for $[v]\in{\rm Comb}(M,\pp)$,
$\varphi:\pi_1(X)\to\Lambda_*$, $\hbasis$ and $\orient$ as usual, a torsion
$\tau_0^\varphi([v],\hbasis,\orient)$ is well-defined, without taking
further quotients.}\end{rem}

As already remarked, $\tau_{{\rm w},0}^\varphi(P,\hbasis,\orient)$ includes the
cases considered by Turaev in~\cite{turaev:Euler} and~\cite{turaev:spinc},
namely manifolds  which are closed or bounded by tori with field pointing
inwards. However, it is not completely obvious that the invariants are exactly
the same. We will now explain more carefully the relation of our construction
with Turaev's work. For the sake of simplicity, we confine ourselves to the
closed case.

Turaev's definition of torsion goes as follows. He first defines torsions for
combinatorial Euler structures on triangulations, and then he  describes a {\em
universal} procedure to map bijectively the set of combinatorial Euler
structures on any given triangulation of a manifold $M$ onto the set of smooth
Euler structures on $M$ ({\it i.e.}~the set of equivalence classes of 
non-singular vector fields on $M$, under homotopy and local modifications). Our
construction goes in the opposite direction. We start with a vector field and
we use it, together with a spine, to produce a combinatorial Euler structure
(actually, a spider) on the ideal triangulation dual to the spine. 

If we first apply our procedure and then Turaev's one, we obtain by composition
a map from vector fields to smooth Euler structures, and it is not completely
obvious that this map is the canonical projection. The point is that the set of
Euler structures is in a natural way an affine space over $H_1(M;\mz)$, and it
is conceivable that the two approaches lead to different choices of the origin.
Nonetheless, since both constructions are universal and very natural, our guess
is that indeed the map described is the projection, or minus it.

We view Proposition~\ref{Euler:class:prop} as an evidence supporting our guess
that the map
described above is the canonical projection. More precisely, assuming the
map to be the projection, Proposition~\ref{Euler:class:prop} is a formal
consequence of a result in~\cite{turaev:Euler}, according to which the cycle
which appears in the proposition is the Poincar\'e\ dual of the relative
characteristic class of $v$ with respect to $-v$, {\em i.e.}~the
Euler class of $v$. We also remark that Proposition~\ref{Euler:class:prop}
could be used as a basis for a direct proof that also our generalized torsions,
in the acyclic closed case, satisfy the duality property established in
Section~2.7 of~\cite{turaev:spinc}.

Another point deserves some emphasis. As mentioned, Turaev's
torsions are (by construction) invariant under local modifications of the field.
It is conceivable that the same is true for our generalized torsions, but
we believe that a direct proof, using only spines, may be hard.

\subsection{Torsion of a (pseudo-)Legendrian link} We will denote by $\leg$ the
set of all equivalence classes of triples $(M,v,L)$ where $M$ is closed, $v$ is
a field on $M$, and $L$ is a link embedded in $M$ and transversal to $v$.
The equivalence relation is generated by homeomorphisms of manifolds and
homotopy of $v$ through fields transverse to $L$. Note that if one allows $L$ to
move isotopically during the homotopy of $v$, the same equivalence relation is
defined. 

For $(M,v,L)$ as above, we will call $L$ a {\em pseudo-Legendrian} link in
$(M,v)$. Our terminology is due to the following example, which also serves as
the main motivation for the definition. If $\xi$ is a cooriented contact
structure on $M$ and $L$ is a Legendrian link in $(M,\xi)$, viewed up to Legendrian
isotopy, then an element $[M,\xi^+,L]\in\leg$ is well-defined, where $\xi^+$ is any
field positively transversal to $\xi$. 

\begin{lem}\label{dig:well-def} A map $D:\leg\to\combd$ can be defined as
follows: given $[M,v,L]$, consider a small open regular neighbourhood $N$
of $L$, and set $D([M,v,L])=[M\setminus N,v]$.
\end{lem}

\dim{dig:well-def} If $N$ is small, $v$ has exactly two concave tangency
circles on each boundary component of $\partial(M\setminus N)$. Independence of
$N$ is immediate, and independence of the representative of $[M,v,L]$ follows
from the fact that any homotopy of $v$ through fields transversal to $L$ can be
replaced by a homotopy which is constant near $L$.\finedim{dig:well-def}

This lemma implies that for $[M,v,L]\in\leg$ we can define torsion invariants
$\tau_0^\varphi(P,\hbasis,\orient)$ for any branched spine $P$ of $D([M,v,L])$.
As usual, one should take into account the action of a mapping class group (but
see also Subsection~\ref{Teich:version:subsection}).

\begin{rem}\label{D:non-injective}{\em The map $D:\leg\to\combd$ defined in the
lemma is neither surjective nor injective. The image of $D$ consists precisely
of pairs $[B,w]$ such that $\partial B$ consists of tori, and each torus is
split into a white and a black annulus. The reason for non-injectivity of $D$
is that several non-isomorphic combed Dehn fillings on such a manifold $[B,w]$ can be
compatible with the black-white splitting.}\end{rem}

\begin{rem}{\em Instead of considering pseudo-Legendrian links in closed combed
manifolds, we may have taken them in manifolds with non-empty boundary and
concave combings. All definitions and constructions are easily
adapted.}\end{rem}

\begin{rem}{\em Since every (pseudo-)Legendrian link has a natural framing,
there is a map which associates to an element of $\leg$ a pair consisting of a
combing and an isotopy class of framed link on the same manifold. One easily
sees that this map is not injective. However, it is surjective. For instance,
given a combing and a framed link, one can realize the link as a
Legendrian one in the unique overtwisted contact structure transversal to the
combing. Note that, by the work of Bennequin~\cite{bennequin}, one cannot use
tight structures. }\end{rem}

\begin{rem}{\em The flexibility of overtwisted contact structures suggests there
could be essentially no difference between pseudo-Legendrian links in a combed
3-manifold and Legendrian links in the unique overtwisted contact structure
homotopic to the combing. To be precise, it is maybe possible use the techniques 
of Eliashberg~\cite{eliash} to answer in the positive to the
following question: {\em for $j=0,1$ let $\xi_j$ be an
overtwisted contact structure on $M$, and let $L_j$ be Legendrian with respect
to $\xi_j$; if $[M,\xi_0^+,L_0]=[M,\xi_1^+,L_1]$ in $\leg$, does there exists an
isotopy of $M$ which carries $\xi_0$ to $\xi_1$ and $L_0$ to $L_1$?}
The content of~\cite{eliash} is a positive answer
for empty links.}\end{rem}

\begin{rem}\label{invariant:of:exterior}{\em  As it often happens with link
invariants, our torsions are actually invariants of the {\em complement} of a
Legendrian link, rather than the Legendrian embedding itself (compare also with
the non-injectivity of $D$ already pointed out above). Therefore, it is clear
that there will be many inequivalent elements of $\leg$ with the same
invariants. In general, the typical situation in which one uses invariants of
exteriors is to compare links in the same manifold. In our setting, we can for
instance use torsions to distinguish different Legendrian links in one and the
same combed manifold. A refinement of this situation will be discussed in the
next subsection.}\end{rem}

\subsection{An embedding-refined of torsion for
links}\label{Teich:version:subsection} As already mentioned, one can take a
Teichm\"uller-type and a moduli-type approach to torsions. We will describe now
a situation where a canonical marking arises, so that the former approach is
more natural, {\it i.e.}~the action on torsions of the mapping class group can
be neglected. 

Let us fix a definite closed manifold $M$. Let $v_0$ and $v_1$ be fields on
$M$, and consider links $L_0$ and $L_1$ transversal to $v_0$ and $v_1$
respectively, so that framings are naturally defined on $L_0$ and $L_1$. We
will assume that $L_0$ and $L_1$ are isotopic in $M$, and we will define
torsions which potentially can prove that they are
{\em not} isotopic {\em as framed links}. 

Let $E(L_i)$ be the complement in $M$ of an open regular neighbourhood of
$L_i$, chosen so that $v_i$ behaves on each torus boundary in the usual way
(with two parallel non-trivial tangency lines). Consider a branched spine $P_i$
which represents $(E(L_i),v_i|_{E(L_i)})$ in the sense of
Theorem~\ref{main:improvement:teo}, and remark that $X(P_i)$ is homeomorphic to
$E(L_i)/(\partial E(L_i))$. By assumption, there is a homeomorphism $f:E(L_0)\to
E(L_1)$ which is the restriction of an automorphism of $M$ isotopic to the
identity. Such an $f$ defines a homeomorphism $F:X(P_0)\to X(P_1)$.

\begin{prop}\label{Teich:version:prop}
If there exists a continuous family $\{v_t\}_{t\in[0,1]}$ of fields on $M$ and
an isotopy of links $\{L_t\}_{t\in[0,1]}$ with $L_t$ transversal to $v_t$ for
all $t$, then
for all choices of $\varphi,\hbasis,\orient$ for $X(P_0)$ we have:
$$\tau_0^\varphi(P_0,\hbasis,\orient)=
\tau_0^{F_*(\varphi)}(P_1,F_*(\hbasis),F_*(\orient)).$$
\end{prop}

\dim{Teich:version:prop} Let $\{Q_i\}_{i=0}^k$ be a sequence of branched
spines, such that $Q_0\cong P_0$, $Q_k\cong P_1$, and $Q_{i-1}\move{\ } Q_i$ is
a sliding move. This move can actually be realized in $M$, so $X(Q_{i-1})$ and
$X(Q_i)$ can be identified using the quotient of  an isotopy of the
complements. The composite homeomorphism $X(P_0)\to X(P_1)$ is isotopic to $F$,
and the conclusion follows from Theorem~\ref{sliding:invariance} and the
accompanying discussion.\finedim{Teich:version:prop}

\subsection{Computational hints}
\paragraph{Spine of a Legendrian link.} As already mentioned in
Remark~\ref{invariant:of:exterior}, the typical setting in which one imagines
to use invariants to distinguish Legendrian links is when a certain closed
combed manifold $(M,v)$ is given, and one restricts to links in
$(M,v)$. In this situation, we can also fix a certain branched spine $P$ of
$(M,v)$, and use it to construct spines of link complements. We first note that on $P$ we
can consider link diagrams, requiring crossings to be away from $S(P)$.
Among link diagrams, it is natural to call $\cont^1$ those which never
meet $S(P)$ by going from one sheet to the other sheet of the locally
two-sheeted area. A $\cont^1$ link projection on $P$ naturally defines a
Legendrian link in $(M,v)$. Conversely, since the complement of $P$ in $M$ is
an open ball, one sees that every Legendrian link in $(M,v)$ is represented by
a $\cont^1$ diagram on $P$. With a little more effort the following can be
established.

\begin{lem}
Every Legendrian link in $(M,v)$ is represented by
a $\cont^1$ diagram on $P$ without crossings.
\end{lem}

\begin{figure}
\centerline{\psfig{file=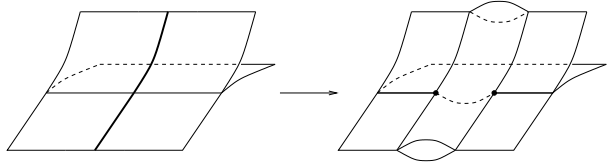,width=8cm}}
\caption{\label{dig:tunnel} How to dig a tunnel in a spine.}
\end{figure}

Now, the tunnel-digging process which takes from a Legendrian link to an element
of $\combd$ can be easily carried over at the level of spines, as suggested in
Fig.~\ref{dig:tunnel}.
This is particularly easy when there are no crossings, but it works in general.
Note that the spine which results from the digging may occasionally be non-standard,
but it is standard as soon as the diagram is complicated enough ({\em e.g.}~if
there are both a crossing and an intersection with $S(P)$).

\paragraph{Boundary operators.} The key point for the computation of torsion is
the knowledge of the boundary operators in the twisted chain complex
$\Cphi_*(P)$. The first step to determine these operators is to describe the
universal cover of $X(P)$, or the maximal Abelian cover when $\Lambda$ is
commutative. Assuming the cover to be perfectly understood, the boundaries in 
$\Cphi_*(P)$ are just liftings of those in $\Ccell_*(X(P);\mz)$. We will now
show that the complex $\Ccell_*(X(P);\mz)$ admits a very easy description,
which seems to indicate that complete calculations should be feasible at least
in some cases, and may be implemented on a computer.

In $X(P)$, we will denote by $\hat{R}$ (respectively $\hat{e}$, $\hat{v}$) the
edge (respectively triangle, tetrahedron) of $\ttt(P)$, dual to a region
(respectively edge, vertex) of $P$. First of all, since there is only one
vertex, we have $\partial\hat{R}=0$ for all $R$. Next, we have
$\partial\hat{e}=\hat{R}_1+\hat{R}_2-\hat{R}_0$, where $R_0,R_1,R_2$ are the
regions incident to $e$, numbered so that $R_1$ and $R_2$ induce on $e$ the same
orientation. (Of course $R_0,R_1,R_2$ need not be different from each other, so
the formula may actually have some cancellation.) Finally,
$\partial\hat{v}=\hat{e}_1+\hat{e}_2-\hat{e}_3-\hat{e}_4$, where $e_1,e_2$ are
the edges which (with respect to the natural orientation) are leaving $v$, and
$e_3,e_4$ are those which are reaching it. (Again, there could be
repetitions.)


\vspace{.3cm}

{\small 

\hspace{8cm} benedett@dm.unipi.it

\hspace{8cm} petronio@dm.unipi.it

\hspace{8cm} Dipartimento di Matematica

\hspace{8cm} Via F.~Buonarroti, 2

\hspace{8cm} I-56127, PISA (Italy)}

\end{document}